\documentclass[moor]{informs3_arxiv}              

\usepackage{natbib}
 \NatBibNumeric
 \def\bibfont{\small}%
 \def\newblock{\ }%
 \bibpunct[, ]{[}{]}{,}{n}{}{,}%

\usepackage[colorlinks=true,breaklinks=true,bookmarks=true,urlcolor=blue,
     citecolor=blue,linkcolor=blue,bookmarksopen=false,draft=false]{hyperref}

\def\EMAIL#1{\href{mailto:#1}{#1}}

\TheoremsNumberedThrough     

\EquationsNumberedThrough    

\usepackage{booktabs} 
\usepackage{xfrac}
\usepackage{algorithm}
\usepackage{algorithmic}
\usepackage[usenames, dvipsnames, table]{xcolor}

\newcommand{\R}{\mathbb{R}}
\newcommand{\Rext}{\mathbb{R}\cup\{+\infty\}}

\newcommand{\abs}[1]{\left\vert#1\right\vert}
\newcommand{\set}[1]{\left\{#1\right\}}
\newcommand{\norm}[1]{\left\Vert#1\right\Vert}

\newcommand{\Eproof}{\hfill $\blacksquare$}
\newcommand{\prox}{\mathrm{prox}}
\newcommand{\dom}[1]{\mathrm{dom}(#1)}
\newcommand{\xb}{x}
\newcommand{\xbar}{\bar{x}}
\newcommand{\yb}{y}
\newcommand{\ub}{u}
\newcommand{\vb}{v}
\newcommand{\cb}{c}
\newcommand{\db}{d}
\newcommand{\Ab}{A}

\newcommand{\Xb}{X}
\newcommand{\Xc}{\mathcal{X}}
\newcommand{\Sc}{\mathbb{S}}

\newcommand{\hb}{h}
\newcommand{\Hb}{H}

\newcommand{\iprods}[1]{\langle #1\rangle}

\newcommand{\intx}[1]{\mathrm{int}\left(#1\right)}

\newcommand{\mblue}[1]{\textcolor{blue}{#1}}

\SingleSpacedXI

\begin{document}


\RUNAUTHOR{Q. Tran-Dinh, A. Kyrillidis, V. Cevher}

\RUNTITLE{A single-phase, proximal path-following framework}

\TITLE{A single-phase, proximal path-following framework\footnote{This is the first revision. The initial version was uploaded on arxiv on March 5, 2016.}}

\ARTICLEAUTHORS{%
\AUTHOR{Quoc Tran-Dinh}
\AFF{Department of Statistics and Operations Research, UNC-Chapel Hill, USA, \EMAIL{quoctd@email.unc.edu}}
\AUTHOR{Anastasios Kyrillidis}
\AFF{University of Texas at Austin, USA, \EMAIL{anastasios@utexas.edu}}
\AUTHOR{Volkan Cevher}
\AFF{Laboratory for Information and Inference Systems (LIONS), EPFL, Switzerland, \EMAIL{volkan.cevher@epfl.ch}}
} 

\ABSTRACT{%
We propose a new proximal, path-following framework for a class of constrained convex problems. We consider settings where the nonlinear---and possibly non-smooth---objective part is endowed with a proximity operator, and the constraint set is equipped with a self-concordant barrier. Our approach relies on the following two main ideas. First, we re-parameterize the optimality condition as an auxiliary problem, such that a good initial point is available; by doing so, a family of alternative paths towards the optimum is generated. Second, we combine the proximal operator with path-following ideas to design a single-phase, proximal, path-following algorithm. 
Our method has several advantages. 
First, it allows handling non-smooth objectives via proximal operators; this avoids lifting the problem dimension in order to accommodate non-smooth components in optimization.
Second, it consists of only a \emph{single phase}: While the overall convergence rate of classical path-following schemes for self-concordant objectives does not suffer from the initialization phase, proximal path-following schemes undergo slow convergence, in order to obtain a good starting point \cite{TranDinh2013e}. 
In this work, we show how to overcome this limitation in the proximal setting and prove that our scheme has the same $\mathcal{O}(\sqrt{\nu}\log(1/\varepsilon))$ worst-case iteration-complexity with standard approaches \cite{Nesterov2004,Nesterov1994} without requiring an initial phase, where $\nu$ is the barrier parameter and $\varepsilon$ is a desired accuracy. 
Finally, our framework allows errors in the calculation of proximal-Newton directions, without sacrificing the worst-case iteration complexity.
We demonstrate the merits of our algorithm via three numerical examples, where proximal operators play a key role.
}%


\MSCCLASS{90C06; 90C25; 90-08}
\ORMSCLASS{Interior-point methods, non-smooth convex programming}

\maketitle

%


%
%
%
\section{Introduction.}\label{sec:intro}
This paper studies the following constrained convex optimization problem: \vspace{-0.15cm}
\begin{equation}\label{eq:constr_cvx}
G^{\star} := \min_{\xb\in\R^p}\Big\{ G(\xb) := \iprods{\cb, \xb} + g(\xb) : \xb \in \Xc \Big\},
\end{equation}   
\vspace{-0.05cm}
where $\cb\in\R^p$, $g$ is a possibly non-smooth, proper, closed and convex function from $\R^p$ to $\R\cup\set{+\infty}$ and $\Xc$ is a nonempty, closed  and convex set in $\R^p$.\footnote{We note that the linear term $\iprods{\cb,\xb}$ can be absorbed into $g$. 
However, we separate it from $g$ for our convenience  in processing numerical examples in the last section.}
We denote by  $\Xc^{\star}$ the optimal solution set of \eqref{eq:constr_cvx}, and by $\xb^{\star}$ an optimal solution in $\Xc^{\star}$.

For convex sets $\Xc$, associated with a self-concordant barrier (see Section \ref{sec:background} for details), and for $G$ \emph{self-concordant and smooth}, \emph{e.g.}, just linear or quadratic, interior point methods (IPMs) often constitute the method-of-choice for solving \eqref{eq:constr_cvx}, with a well-characterized worst-case complexity. 
A non-exhaustive list of instances of \eqref{eq:constr_cvx} includes linear programs, quadratic programs, second order cone programs, semi-definite programs, and geometric optimization  \cite{Andersen2011,Becker2011a,BenTal2001,Gondzio2012,Nemirovski2006,Nemirovski2008,Nesterov2004, Nesterov1997,Renegar2001,Roos2006,vanderbei2013optimization,Wright1997}.

At the heart of IPMs lies the notion of \emph{interior barriers}: these mimic the effect of the constraint set $\Xc$ in \eqref{eq:constr_cvx} by appropriately penalizing  the objective function with a barrier $f$ over the set $\Xc$, as follows:
\begin{equation}{\label{eq:path_following_form}}
F^{\star}_t := \min_{\xb\in\intx{\Xc}}\Big\{ F_t(\xb) :=  \tfrac{1}{t} \cdot G(\xb) +  f(\xb)\Big\}.
\end{equation} 
Here, $f$ models the structure of the feasible set $\Xc$ and $t > 0$ is a penalty parameter.
For different values of $t$, the regularized problem generates a sequence of solutions $\left\{\xb^\star(t) : t > 0 \right\}$, known as the \emph{central path}, converging to $\xb^{\star}$ of \eqref{eq:constr_cvx}, as $t$ goes to $0^{+}$.
\emph{Path-following methods} operate along the central path: 
for a properly decreasing sequence of $t$ values, they solve \eqref{eq:path_following_form} only approximately, by performing a few Newton iterations for each $t$ value;
standard path-following schemes even perform just \textit{one} Newton iteration, assuming a linear objective $G(x) := \langle c, x \rangle$ with no non-smooth term $g(x)$.
For such problem cases, this is sufficient to guarantee that the approximate solution lies sufficiently close to the central path, and operates as warm-start for the next value of $t$ in \eqref{eq:path_following_form} \cite{BenTal2001,nemirovski2001lectures,Nesterov1994,Nesterov2004}.
One requirement is that the initial point must lie within a predefined neighborhood of the central path. 
In their seminal work \cite{Nesterov1994}, Nesterov and Nemirovski showed that such methods admit a polynomial worst-case complexity, as long as the Newton method has polynomial complexity.

Based on the above, standard schemes \cite{Nemirovski2006,Nesterov2004,Nesterov1994}
can be characterized by two phases: \textsc{Phase I} and \textsc{Phase II}. 
In \textsc{Phase I} and for an initial value of $t$, say $t_0$, one has to solve \eqref{eq:path_following_form} carefully in order to determine a good initial point for \textsc{Phase II}; this implies solving \eqref{eq:path_following_form} up to sufficient accuracy, such that the Newton method for \eqref{eq:path_following_form} admits fast convergence. 
In \textsc{Phase II} and using the output of \textsc{Phase I} as a warm-start, we path-follow with a provably polynomial time complexity. 

Taking into account both phases, standard path-following algorithms---where \eqref{eq:path_following_form} is a \emph{self-concordant} objective---are characterized by the following iteration complexity. The total number of iterations required to obtain an $\varepsilon$-solution is
\begin{align}{\label{eq:intro_com}}
\mathcal{O}\left(\sqrt{\nu}\log\left(\tfrac{1}{\varepsilon}\right)\right).
\end{align} 
Here, $\nu$ is a barrier parameter (see Section \ref{sec:background} for details) and $\varepsilon$ is the approximate parameter, according to the following definition:
\begin{definition}\label{de:approx_sol}
Given a tolerance $\varepsilon > 0$, we say that $\xb^{\star}_{\varepsilon}$ is an $\varepsilon$-solution for \eqref{eq:constr_cvx} if
\begin{equation*}
\xb^{\star}_{\varepsilon} \in \Xc, ~~\text{and}~~G(\xb^{\star}_{\varepsilon}) - G^{\star} \leq \varepsilon.
\end{equation*}
\end{definition}

\subsection{Path-following schemes for non-smooth objectives.}
For many applications in machine learning, optimization and signal processing \cite{BenTal2001,Parikh2013,TranDinh2013e}, the $g$ part in \eqref{eq:constr_cvx} could be non-smooth (or even smooth but non-self-concordant).
Such a $g$ term is usually included in the optimization in order to leverage the true underlying structure in $\xb^\star$. 
An example is the $\ell_1$-norm regularization, \emph{i.e.}, $g(\xb) = \|\xb\|_1$, with applications in high-dimensional statistics, compressive sensing, scientific and medical imaging \cite{candes2013phase, gramfort2009improving,jaganathan2013sparse,jenatton2011multiscale,millane1990phase,rapaport2008classification,subramanian2005gene,zhou2010association}, among others. Other examples for $g$ include the indicator function of a convex set \cite{Parikh2013}, the $\ell_{1,2}$-group norm \cite{baldassarre2013group,jenatton2011structured, kyrillidis2015structured}, and the nuclear norm \cite{cai2010singular} using in low-rank matrix approximation.

Unfortunately, non-smoothness in the objective reduces the optimization efficiency.
In such settings, one can often reformulate \eqref{eq:constr_cvx} into a standard conic program, by introducing slack variables and additional constraints to model $g$.
Such a technique is known as \emph{disciplined convex programming} (DCP) \cite{Grant2006} and has been incorporated in well-known software packages, such as CVX \cite{Grant2006} and YALMIP \cite{Loefberg2004}. 
Existing off-the-shelf solvers are then utilized to solve the resulting problem.
However, DCP could potentially increase the problem dimension significantly; this, in sequence, reduces the efficiency of the IPMs.
For instance, in the example above where $g(\xb) = \Vert\xb\Vert_1$, DCP introduces $p$ slack variables to reformulate $g$ into $O(p)$ additional linear constraints; when $g(X) = \|X\|_*$, \emph{i.e.}, $g$ is the nuclear norm (sum of singular values of $X \in \mathbb{R}^{p \times q}$), then it can be \emph{smoothed} via a semi-definite formulation, where the memory requirements and the volume of computation per iteration are high \cite{liu2009interior}.

In this paper, we focus on cases where $g$ is endowed with a generalized proximity operator, associated with a local norm $\norm{\cdot}_{\xb}$ (see Section \ref{sec:background} for details): 
\begin{equation*}
\textrm{prox}_g(\ub) := \argmin_{\vb\in\R^p}\set{ g(\vb) + \sfrac{1}{2} \cdot \Vert \vb - \ub \Vert_{\xb}^2 }.
\end{equation*}
Such proximity operators have been used extensively in non-smooth optimization problems, and proven to be efficient in real applications, under common gradient Lipschitz-continuity and strong convexity assumptions on the objective function \cite{Beck2009, Combettes2011, Nesterov2004}. 
However, for generic $\Xc$ constraints in \eqref{eq:constr_cvx}, the resulting interior barrier $f$ in \eqref{eq:path_following_form} does not have Lipschitz continuous gradients and, thus, prevents us from trivially recycling such ideas.
This necessitates the design of a new class of path-following schemes, that exploit proximal operators and thus can accommodate non-smooth terms in the objective.

To the best of our knowledge, \cite{TranDinh2013e} is the first work that treats jointly interior barrier path-following schemes and proximity operators, in order to construct new proximal path-following algorithms for problems as in \eqref{eq:constr_cvx}. 
According to \cite{TranDinh2013e}, the proposed algorithm follows a two-phase approach, with \textsc{Phase II} having the same worst-case iteration-complexity as in \eqref{eq:intro_com} (up to constants) \cite{Nesterov2004,Nesterov1994}. 
However, the initialization \textsc{Phase I} in \cite{TranDinh2013e} requires substantial computational effort, which usually dominates the overall computational time. 
In particular, to find a good initial point, \cite{TranDinh2013e} uses a damped-step proximal-Newton scheme for \eqref{eq:path_following_form}, starting from an arbitrary initial point $\xb^0$ and for arbitrary selected $t_0 > 0$. 
For such configuration, \cite{TranDinh2013e} requires
\begin{align*}
\left \lfloor \frac{F_{t_0}(\xb^0) - F_{t_0}(\xb^\star(t_0))}{\omega \left((1 - \kappa)\beta \right)} \right \rfloor
\end{align*} 
damped-step Newton iterations in \textsc{Phase I} in order to find a point close to the optimal solution of \eqref{eq:path_following_form}, say $x^\star(t_0)$, for the selected $t_0$.
Here, $\kappa \in (0,1),~\beta \in (0, 0.15]$, and $\omega(\tau) := \tau - \log(1 + \tau) \geq 0$; see \cite[Theorem 4.4]{TranDinh2013e} for more details. 
\emph{I.e.}, in stark contrast to the global iteration complexity \eqref{eq:intro_com} of smooth path-following schemes, \textsc{Phase I} of \cite{TranDinh2013e} might require a substantial number of iterations just to converge to a point close to the central path, and depends on the arbitrary initial point selection $\xb^0$. 

\subsection{Motivation.}
From our discussion so far, it is clear that most existing works on path-following schemes require two phases. 
In the case of \emph{smooth} self-concordant objectives in \eqref{eq:constr_cvx},
\textsc{Phase I} is often implemented as a damped-step Newton scheme, which has sublinear convergence rate, or an auxiliary path-following scheme, with linear convergence rate that satisfies the global, worst-case complexity in \eqref{eq:intro_com} \cite{Nesterov2004,Nesterov1994}.
In standard conic programming, one can unify a two-phase algorithm in a single-phase IP path-following scheme via  homogeneous and self-dual  embedded strategies; see, \emph{e.g.}, \cite{Sturm1999,Tutunku2003,Wright1997}.
Such strategies parameterize the KKT condition of the primal and dual conic program so that one can immediately have an initial point, without performing \textsc{Phase I}.
So far and to the best of our knowledge, it remains unclear how such an auxiliary path-following scheme can find an initial point for \emph{non-smooth} objectives in \eqref{eq:constr_cvx}. 

\subsection{Our contributions.}
The goal of this paper is to develop a new single-phase, proximal path-following algorithm for \eqref{eq:constr_cvx}. 
To do so, we first re-parameterize the optimality condition of the barrier problem associated with \eqref{eq:constr_cvx} as a \emph{parametric monotone inclusion} (PMI).
Then, we design a proximal path-following scheme to approximate the solution of such PMI, while controlling the penalty parameter. 
Finally, we show how to recover an approximate solution of \eqref{eq:constr_cvx}, from the approximate solution of the PMI. 

The main contributions of this paper can be summarized as follows:
\begin{itemize}
\item[$(i)$] We introduce a new parameterization for the optimality condition of \eqref{eq:path_following_form} to appropriately select the parameters such that less computation for initialization is needed. 
Thus, with an appropriate choice of parameters, we show how we can eliminate the slowly-convergent  \textsc{Phase I} in \cite{TranDinh2013e}, while we still maintain the global, polynomial time, worst-case iteration-complexity.

In particular, we propose novel---checkable \emph{a priori}---conditions over the set of initial points that can achieve the desiderata; this, in turn, provides rigorous configurations of the algorithm's parameters such that the worst-case iteration complexity guarantee is obtained provably, avoiding the slowly convergent initialization procedures proposed so far for non-smooth optimization in \eqref{eq:constr_cvx}. 

\item[$(ii)$] We design a single-phase, path-following algorithm to compute an $\varepsilon$-solution of \eqref{eq:constr_cvx}. For each $t$ value, the resulting algorithm only requires a \emph{single approximate Newton iteration} (see \cite{TranDinh2013e}), followed by a proximal step, of a strongly convex quadratic composite subproblem.
We will use the term \emph{proximal Newton step} when referring to these two steps.
The algorithm allows inexact Newton steps, with a verifiable stopping criterion (\emph{cf.} eq. \eqref{eq:approx_subsol2}). 
\end{itemize}

In particular, we establish the following result:
\begin{theorem}
The total number of proximal Newton iterations required, in order to reach an $\varepsilon$-solution of \eqref{eq:constr_cvx}, 
is upper bounded by
$\mathcal{O}\left(\sqrt{\nu}\log\left(\frac{\nu}{\varepsilon}\right)\right)$. 
\end{theorem} 
A complete and formal description of the above theorem and its proof are provided in Section \ref{sec:proxPF_alg}. 
Our \emph{proximal} algorithm admits the same iteration-complexity, as standard path-following methods \cite{Nesterov2004,Nesterov1994} (up to a constant). 
To highlight the iteration complexity gains from the two-phase algorithm in \cite[Theorem 4.4]{TranDinh2013e}, recall that in the latter case, the total number of proximal Newton steps are bounded by:
\begin{align*}
\left \lfloor \frac{F_{t_0}(\xb^0) - F_{t_0}(\xb^\star_{t_0})}{\omega \left((1 - \kappa)\beta \right)} \right \rfloor + \mathcal{O}\left(\sqrt{\nu}\log\left(\frac{\nu}{\varepsilon}\right)\right),
\end{align*} 
where the first term is in \textsc{Phase I} as mentioned previously, and the second one is  in \textsc{Phase II}.

Our algorithm requires a well-chosen initial point that avoids \textsc{Phase I}; one such case is that of an approximation of the analytical center $\bar{\xb}^{\star}_f$ of the barrier $f$ (see Section \ref{sec:background} for details). 
In the text, we argue that evaluating this point is much easier than finding an initial point $\xb^0$ using \textsc{Phase I}, as in \cite{TranDinh2013e}.
In addition, for many feasible sets $\Xc$ in \eqref{eq:constr_cvx}, we can explicitly and easily compute $\bar{\xb}^{\star}_f$ of $f$ (see Section~\ref{sec:num_exams} for examples).

\subsection{The structure of the paper.}
This paper is organized as follows.
Sections \ref{sec:background} and \ref{sec:central_path} contain basic definitions and notions, used in our analysis. 
We introduce  a new re-parameterization of the central path in order to obtain a \textit{predefined} initial point.
Section \ref{sec:proxPF_alg} presents a novel algorithm and its complexity theory for the non-smooth  objective function.
Section \ref{sec:num_exams} provides three numerical examples that highlight the merits of our algorithm. 

\section{Preliminaries.}\label{sec:background}
In this section, we provide the basic notation used in the rest of the paper, as well as two key concepts: proximity operators and self-concordant (barrier) functions.

\subsection{Basic definitions.} 
Given $\xb, \yb \in \mathbb{R}^p$, we use $\iprods{\xb, \yb}$ or $\xb^T\yb$ to denote the inner product in $\R^p$. 
For a proper, closed and convex function $g$, we denote by $\dom{g}$ its domain, (\emph{i.e.}, $\dom{g} := \set{\xb\in\R^p : g(\xb) < + \infty})$,
and by $\partial{g}(\xb) := \set{\vb\in\R^p | g(\yb) \geq g(\xb) + \iprods{\vb, \yb - \xb}, ~\forall\yb\in\dom{g} }$ its subdifferential at $\xb$. 
We also denote by $\mathrm{Dom}(g) := \mathrm{cl}(\dom{g})$ the closure of $\dom{g}$ \cite{Rockafellar1970}.
We use $\mathcal{C}^3(\mathcal{X})$ to denote the class of three times continuously differentiable functions from $\mathcal{X}\subseteq\mathbb{R}^p$ to $\mathbb{R}$.

For a given twice differentiable function $f$ such that $\nabla^2f(\xb) \succ 0 $ at some $\xb\in\dom{f}$, we define the local norm, and its dual, as 
\begin{equation*}
\norm{\ub}_{\xb} := \iprods{\nabla^2 f(\xb)\ub, \ub}^{1/2}, \forall \ub\in\R^p, ~~\text{and}~~\norm{\vb}_{\xb}^{*} := \max_{\norm{\ub}_{\xb} \leq 1}\iprods{\ub,\vb} = \iprods{\nabla^2f(\xb)^{-1}\vb, \vb}^{1/2},
\end{equation*} respectively, 
for $\ub, \vb\in\mathbb{R}^p$. 
Note that the Cauchy-Schwarz inequality holds, \emph{i.e.}, $\iprods{\ub, \vb} \leq \norm{\ub}_{\xb}\norm{\vb}^{*}_{\xb}$.

\subsection{Generalized proximity operators.}
The generalized proximity operator of a proper, closed and convex function $g$ is defined as the following program:
\begin{equation}\label{eq:prox_oper}
\textrm{prox}_g(\ub) := \argmin_{\vb\in\R^p}\set{ g(\vb) + \sfrac{1}{2} \cdot \Vert \vb - \ub \Vert_{\xb}^2 }.
\end{equation}
When $\nabla^2{f} = \mathbb{I}$---the identity matrix---in the local norm, \eqref{eq:prox_oper} becomes a standard proximal operator \cite{Bauschke2011}. 
Computing $\textrm{prox}_g$ might be hard even for such cases. 
Nevertheless, there exist structured smooth and non-smooth convex functions $g$ with $\textrm{prox}_g$ that comes with a closed-form solution or can be computed with low computational complexity.
We capture this idea in the following definition.

\begin{definition}[\textit{Tractable proximity operator}]\label{de:tractable_proximity}
A proper, closed and convex function $g : \R^p\to\Rext$ has a \emph{tractable} proximity operator if \eqref{eq:prox_oper} can be computed efficiently via a closed-form solution or via a polynomial time algorithm.
\end{definition}

Examples of such functions include the $\ell_1$-norm---where the proximity operator is the well-known soft-thresholding operator \cite{Combettes2011}---and the indicator functions of simple sets (\emph{e.g.}, boxes, cones and simplexes)---where the proximity operator is simply the projection operator.
Further examples can be found in \cite{Bauschke2011,Combettes2011,Parikh2013}.
Observe that, due to the existence of a closed-form solution for most well-known proximity operators, one can always compute $\textrm{prox}_{\tfrac{1}{t}\cdot g}$ efficiently and its computational complexity does not depend on the value of the regularization parameter $t$. 
Our main result does not require the tractability of computing the proximity operator of $g$; 
it will be used to analyze the overall computational complexity in Subsection~\ref{subsec:overall_comp_complexity}.

Some properties of the proximity operator $\textrm{prox}_g$ are described in the next lemma:
\begin{lemma}{\label{lem:proximal_prop}}
The generalized proximal operator defined in \eqref{eq:prox_oper} is co-coercive and therefore nonexpansive w.r.t. the local norms, \emph{i.e.},
\begin{align}
\mathrm{[co\textrm{-}coercive]:} ~~~~~& \left({\rm prox}_g(\ub) - {\rm prox}_g(\vb)\right)^T(\ub - \vb) \geq \norm{{\rm prox}_g(\ub) - {\rm prox}_g(\vb)}_{\xb}^2, \label{eq:cocoercive}\\
\mathrm{[nonexpansive]:}~~ &\norm{{\rm prox}_g(\ub) - {\rm prox}_g(\vb)}_{\xb} \leq \norm{\ub - \vb}_{\xb}^{*}, ~\forall \ub, \vb \in \mathbb{R}^p. \label{eq:nonexapansiveness}
\end{align}
\end{lemma} 
The proof to this lemma can be found in \cite[Lemma~2]{Tran-Dinh2013a}.

\subsection{Self-concordant functions and self-concordant barriers.}
A concept used in our analysis is the self-concordance property, introduced by Nesterov and Nemirovskii \cite{Nesterov2004,Nesterov1994}.

\begin{definition}\label{de:concordant}
A univariate convex function $\varphi \in \mathcal{C}^3(\dom{\varphi})$ is called \emph{standard self-concordant} if 
$\abs{\varphi'''(\tau)} \leq 2\varphi''(\tau)^{3/2}$  for all $\tau\in\dom{\varphi}$, where $\dom{\varphi}$ is an open set in $\R$. 
Moreover, a function $f:\dom{f}\subseteq \R^p \to \R$ is standard self-concordant if, for any $\xb \in\dom{f}$ and $\vb \in\R^p$, the univariate function $\varphi$ where $\tau \mapsto \varphi(\tau) := f(\xb + \tau \vb)$ is standard self-concordant.
\end{definition}

\begin{definition}\label{de:self_con_barrier}
A standard self-concordant function $f : \dom{f}\subset\R^p\to\R$ is a \emph{$\nu$-self-concordant barrier} for the set $\mathrm{Dom}(f)$ with parameter $\nu > 0$, if 
\begin{equation*}
\sup_{\ub \in\R^p} \left \{2\iprods{\nabla{f}(\xb),\ub} - \Vert \ub\Vert_{\xb}^2\right \} \leq \nu, ~~\forall \xb\in\dom{f}.
\end{equation*} 
In addition, $f(x)\to\infty$ as $x$ tends to the boundary of $\dom{f}$.
\end{definition}

We note that when $\nabla^2 f$ is non-degenerate (particularly, when  $\dom{f}$ contains no straight line \cite[Theorem 4.1.3.]{Nesterov2004}), a $\nu$-self-concordant function $f$ satisfies
\begin{align}{\label{eq:used}}
\norm{\nabla f(\xb)}_{\xb}^{\ast} \leq \sqrt{\nu}, ~~\forall \xb\in\dom{f}.
\end{align}
Self-concordant functions have non-global Lipschitz gradients and can be used to analyze the complexity of Newton-methods \cite{Boyd2004,Nesterov2004,Nesterov1994}, as well as first-order variants \cite{frostig2014competing}.
For more details on self-concordant functions and self-concordant barriers, we refer the reader to Chapter 4 of \cite{Nesterov2004}. 

Several simple sets are equipped with a self-concordant barrier. 
For instance, $f_{\R^p_{+}}(\xb) := -\sum_{i=1}^n\log(x_i)$ is an $n$-self-concordant barrier of the orthant cone $\R^p_{+}$, 
$f(\xb, t) = -\log(t^2 - \norm{\xb}_2^2)$ is a $2$-self-concordant barrier of the Lorentz cone $\mathcal{L}_{n+1} := \set{(\xb, t) \in\R^p\times\mathbb{R}_{+} : \norm{\xb}_2 \leq t}$, 
and the semidefinite cone $\mathbb{S}^n_{+}$ is endowed with the $n$-self-concordant barrier $f_{\mathbb{S}_{+}^n}(\Xb) := -\log\det(\Xb)$. 
In addition, other convex sets, such as hyperbolic polynomials and convex cones, are also characterized by explicit self-concordant barriers \cite{lewis2001self,Nesterov2006f}.
Generally, any closed and convex set---with nonempty interior and not containing a straight line---is endowed with a self-concordant barrier; see \cite{Nesterov2004,Nesterov1994}. 

Finally, we define the analytical center $\xbar^{\star}_f$ of $f$ as 
\begin{equation}\label{eq:analytical_center}
\xbar^{\star}_f := \argmin\set{ f(\xb) : \xb\in\intx{\Xc} }~~\Leftrightarrow~~\nabla{f}(\xbar^{\star}_f) = 0.
\end{equation}
If $\Xc$ is bounded, then $\xbar^{\star}_f$ exists and is unique \cite{Nesterov2011c}. 
Some properties of the analytical center are presented in Section \ref{sec:central_path}. 
In this paper, we develop algorithms for \eqref{eq:constr_cvx} with a general self-concordant barrier $f$ of $\Xc$ as defined by Definition \ref{de:self_con_barrier}.

\subsection{Basic assumptions.} 
We make the following assumption, regarding problem \eqref{eq:constr_cvx}. 
\begin{assumption}\label{as:A1}
The solution set $\Xc^{\star}$ of \eqref{eq:constr_cvx} is nonempty.
The objective function $g$ in \eqref{eq:constr_cvx} is proper, closed and convex, and $\Xc\cap\dom{g}\neq\emptyset$.
The feasible set $\Xc$ is nonempty, closed and convex $($with nonempty interior $\intx{\Xc}$$)$ and is endowed with a $\nu$-self-concordant barrier $f$ such that $\mathrm{Dom}(f) := \textrm{cl}(\dom{f}) = \Xc$.
The analytical center $\xbar^{\star}_f$ of $f$ exists.
\end{assumption}

Except for the last condition, Assumption \ref{as:A1} is common for interior-point methods. 
The last condition can be satisfied by adding an auxiliary constraint $\norm{\xb}_2 \leq R$ for sufficiently large $R$; this technique has been also used in \cite{Nesterov1994} and it does not affect the solution of \eqref{eq:constr_cvx} when $R$ is large.

\section{Re-parameterizing the central path.}\label{sec:central_path}
In this section, we introduce a new parameterization strategy, which will be used in our scheme for \eqref{eq:constr_cvx}. 

\subsection{Barrier formulation and central path of \eqref{eq:constr_cvx}.}
Since $\Xc$ is endowed with a $\nu$-self-concordant barrier $f$, according to Assumption A.\ref{as:A1}, the barrier formulation of \eqref{eq:constr_cvx} is given by
\begin{equation}\label{eq:barrier_prob2}
F^{\star}_t := \min_{\xb\in\intx{\Xc}}\Big\{ F_t(\xb) := \tfrac{1}{t} \cdot G(\xb) + f(\xb) \equiv \tfrac{1}{t} \cdot \left(\iprods{\cb,\xb} + g(\xb)\right) + f(\xb) \Big\},
\end{equation}
where $t > 0$ is the penalty parameter.
We denote by $\xbar^{\star}_t$ the solution of \eqref{eq:barrier_prob2} at a given value $t > 0$. 
Define $r_{t}(\xb) := \tfrac{1}{t} \cdot \left(\cb + \partial{g}(\xb)\right) +  \nabla f (\xb)$.
The optimality condition of \eqref{eq:barrier_prob2} is necessary and sufficient for $\xbar^{\star}_t$ to be an optimal solution of \eqref{eq:barrier_prob2}, and can be written as follows:
\begin{equation}\label{eq:opt_cond2}
0 \in  r_t(\xbar^\star_t) \equiv \tfrac{1}{t} \cdot \left(\cb + \partial{g}(\bar{\xb}^{\star}_t)\right) + \nabla{f}(\bar{\xb}^{\star}_t).
\end{equation}
We also denote by $\bar{\mathcal{C}} := \set{ \xbar^{\star}_t : t > 0}$ the set of solutions of \eqref{eq:barrier_prob2}, which generates a central path (or a solution trajectory) associated with \eqref{eq:constr_cvx}. 
We refer to each solution $\xbar^{\star}_t$ as a central point.

\subsection{Parameterization of the optimality condition.} 
Let us fix $\xb^0\in\dom{f}$; a specific selection of $\xb^0$ is provided later on. 
For given $\xb^0$, let $\xi_0 \in \partial{g}(\xb^0)$ be an arbitrary subgradient of $g$ at $\xb^0$, and set $\zeta_0 := \nabla{f}(\xb^0) + \tfrac{1}{t_0} \cdot (\cb + \xi_0)$. 
For a given parameter $\eta > 0$, define
\begin{align}\label{eq:rteta}
h_{\eta}(\xb) := f(\xb) - \eta\iprods{\zeta_0, \xb} \quad \text{and} \quad r_{t, \eta}(\xb) := \tfrac{1}{t} \cdot \left(\cb + \partial{g}(\xb)\right) +  \nabla{h}_{\eta}(\xb).
\end{align} 
with the gradient $\nabla{h}_{\eta}(\xb) := \nabla{f}(\xb) - \eta\zeta_0$. 
We further define an $\eta$-parameterized version of \eqref{eq:barrier_prob2} as
\begin{equation}\label{eq:modified_cvx_prob}
H^{\star}_t := \min_{\xb\in\intx{\Xc}}\set{ H_t(\xb) := \tfrac{1}{t} \cdot \left(\iprods{\cb, \xb} + g(\xb)\right) + h_{\eta}(\xb)}. 
\end{equation} 
We denote by $\xb^{\star}_t$ the solution of \eqref{eq:modified_cvx_prob}, given $t > 0$. 
Observe that, for a fixed value of $\eta > 0$, the optimality condition of \eqref{eq:modified_cvx_prob} at $\xb^\star_t$ is given by
\begin{equation}\label{eq:param_opt_cond2}
0 \in r_{t, \eta}(\xb^{\star}_t) \equiv \tfrac{1}{t} \cdot \left(\cb + \partial{g}(\xb^{\star}_t)\right) + \nabla{h}_{\eta}(\xb^{\star}_t).
\end{equation} 

\medskip
Next, we provide some remarks regarding the $\eta$-parameterized problem in \eqref{eq:modified_cvx_prob}:
\begin{enumerate}
\item Clearly, if we set $\eta = 0$, $h_\eta(\xb) \equiv f(\xb)$ and thus, \eqref{eq:modified_cvx_prob} is equivalent to \eqref{eq:barrier_prob2}. Therefore, for any other value $\eta > 0$, the problem in \eqref{eq:modified_cvx_prob} differs from the original formulation \eqref{eq:barrier_prob2} by a factor $- \eta \langle \zeta_0, \xb\rangle$.
\item Fix parameters $\eta > 0, t > 0$ and let $\xb^{\star}_t$ be the solution of \eqref{eq:modified_cvx_prob}, which is different from the solution $\bar{\xb}^{\star}_t$ of \eqref{eq:barrier_prob2}, given the remark above. However, as $t \rightarrow 0$ in a path-following scheme, both $\xb^{\star}_t$ and $\bar{\xb}^{\star}_t$ converge to an optimum $\xb^\star$ of \eqref{eq:constr_cvx}. 
\item Based on the above, for fixed $t > 0$ and different values of $\eta$, \eqref{eq:modified_cvx_prob} leads to a family of paths towards $\xb^\star$ of \eqref{eq:constr_cvx}. 
\end{enumerate}

Our aim in this paper is to properly combine the quantities $t_0$, $\xb^0$ and $\eta$, such that 
$(i)$ solving iteratively \eqref{eq:modified_cvx_prob} always has fast convergence (even at the initial point $\xb^0$) and, 
$(ii)$ while \eqref{eq:modified_cvx_prob} differs from \eqref{eq:barrier_prob2}, its solution trajectory is closely related to the solution trajectory of the original barrier formulation. The above are further discussed in the next subsections.

\subsection{A functional connection between solutions of \eqref{eq:barrier_prob2} and \eqref{eq:modified_cvx_prob}.}
Given the definitions above, let us first study the relationship between \emph{exact} solutions of \eqref{eq:barrier_prob2} and \eqref{eq:modified_cvx_prob}, for fixed values $t > 0$ and $\eta > 0$. 
\begin{lemma}\label{le:Delta_star_est}
Let $t > 0$ be fixed. Assume $\eta > 0$ and $\zeta_0$ is chosen such that $m_0 = \eta \Vert\zeta_0 \Vert_{\xbar^{\star}_t}^{*} < 1$. Define $\bar{\Delta}_t := \Vert x^{\star}_t - \bar{x}^{\star}_t\Vert_{\bar{x}^{\star}_t}$ as the local distance between $\bar{x}^{\star}_t$ and $x^{\star}_t$, the solutions of  \eqref{eq:barrier_prob2} and \eqref{eq:modified_cvx_prob}, respectively. Then,
\begin{equation*}
\bar{\Delta}_t \leq \frac{m_0}{1 - m_0}.
\end{equation*}
\end{lemma}

\textit{Proof.} 
Let $x^{\star}_{t}$ be the solution of \eqref{eq:modified_cvx_prob} and $\bar{x}^{\star}_t$ be the solution of \eqref{eq:barrier_prob2}. 
By the optimality conditions in \eqref{eq:opt_cond2} and \eqref{eq:param_opt_cond2}, we have $-t\nabla{f}(\bar{x}^{\star}_t) \in \partial{G}(\bar{x}^{\star}_t)$ and $-t\nabla{h_{\eta}}(x^{\star}_t) \in \partial{G}(x^{\star}_t)$.
Moreover, by the convexity of $G$, we have $\iprods{\nabla{f}(\bar{x}^{\star}_t) - \nabla{h_{\eta}}(x^{\star}_t), x^{\star}_t - \bar{x}^{\star}_t} \geq 0$. 
Using the definition $\nabla{h}_{\eta}(\xb) := \nabla{f}(\xb) - \eta\zeta_0$, the last inequality  leads to
\begin{align*}
\iprods{\nabla{f}(x^{\star}_t) - \nabla{f}(\bar{x}^{\star}_t), x^{\star}_t - \bar{x}^{\star}_t} \leq \eta\iprods{\zeta_0, x^{\star}_t - \bar{x}^{\star}_t}.
\end{align*}
Further, by \cite[Theorem 4.1.5]{Nesterov2004} and the Cauchy-Schwarz inequality, this inequality implies
\begin{equation}\label{eq:lm52_est1}
\frac{\Vert x^{\star}_t - \bar{x}^{\star}_t\Vert_{\bar{x}^{\star}_t}}{1 +  \Vert x^{\star}_t - \bar{x}^{\star}_t\Vert_{\bar{x}^{\star}_t}} 
\leq  \eta\Vert \zeta_0 \Vert_{\bar{x}^{\star}_t}^{*} \quad \Longrightarrow \quad \bar{\Delta}_t \leq \frac{m_0}{1-m_0}, 
\end{equation} which completes the proof of this lemma. 
\Eproof

\medskip
The above lemma indicates that, while \eqref{eq:barrier_prob2} and \eqref{eq:modified_cvx_prob} define different central paths towards $\xb^\star$, there is an upper bound $m_0$ on the distance between  $\bar{x}^{\star}_t$ and $x^{\star}_t$, which is controlled by the selection of $\eta, t_0$ and $x^0$. 
However, $\Vert\zeta_0 \Vert_{\xbar^{\star}_t}^{*}$ cannot be evaluated a priori, since $\xbar^{\star}_t$ is unknown. 

\subsection{Estimate upper bound for $\bar{\Delta}_t$.}
We can overcome this difficulty by using an approximation of the \emph{analytical center point} $\bar{\xb}_f^\star$ in \eqref{eq:analytical_center}. 
A key property of $\xbar^\star_f$ is the following \cite[Corollary 4.2.1]{Nesterov2004}: 
Define $n_{\nu} := \nu+2\sqrt{\nu}$, where $\nu$ is the self-concordant barrier parameter. 
If $f$ is a logarithmically homogeneous self-concordant barrier, then we set $n_{\nu} := 1$ \cite{Nesterov2004}.
Then, $\Vert \vb\Vert_{\xb}^{*} \leq n_{\nu} \Vert\vb\Vert_{\xbar_f^{\star}}^{*}$ for any $\xb\in\mathrm{int}(\Xc)$ and $\vb\in\R^p$. 
This observation leads to the following Corollary; the proof easily follows from that of Lemma \ref{le:Delta_star_est} and the properties above.

\begin{corollary}{\label{co:Delta_star_est}}
Consider the configuration in Lemma \ref{le:Delta_star_est} and define $\bar{m}_0 = \eta n_\nu \Vert\zeta_0 \Vert_{\xbar^{\star}_f}^{*} < 1$. Then, 
\begin{equation}\label{eq:Delta_star_est}
\bar{\Delta}_t \leq \frac{\bar{m}_0}{1 - \bar{m}_0}, \quad \quad \forall t > 0.
\end{equation}
Moreover, if we choose the initial point $\xb^0$ as $\xb^0 := \xbar^{\star}_f$, then $\bar{m}_0 = \tfrac{n_{\nu} \eta}{t_0} \cdot \Vert c + \xi_0\Vert_{\xbar^{\star}_f}^{*} \geq m_0$, where $m_0$ is defined in Lemma \ref{le:Delta_star_est}.
\end{corollary}

\textit{Proof.} 
By \cite[Corollary 4.2.1]{Nesterov2004}, one observes that $\Vert \zeta_0 \Vert_{\bar{x}^{\star}_t}^{*} \leq n_{\nu} \Vert \zeta_0 \Vert_{\bar{x}^{\star}_f}^{*}$, where $\xbar^{\star}_f$ is the analytical center of $f$. 
Following the same motions with the proof of Lemma \ref{le:Delta_star_est}, we obtain \eqref{eq:Delta_star_est}. 
Further, using the property $\nabla f(\xbar^{\star}_f) = 0$ and the definition of $\zeta_0$, we obtain the last statement.
\Eproof

\medskip
In the corollary above, we bound the quantity $\bar{m}_0$ using the local norm at the analytical center $\bar{\xb}^{\star}_f$.
This will allow us to estimate the theoretical worst-case bound in Theorem \ref{th:worst_case_complexity2}, described next. 
This corollary also suggests to choose an initial point $\xb^0 := \bar{\xb}^{\star}_f$, assuming  $\bar{\xb}^{\star}_f$ is available or easy to compute \emph{exactly}. 
Later in the text, we propose initialization conditions that can be checked \emph{a priori} and, if they are satisfied, such initial points are sufficient for our results to apply. \emph{E.g.}, consider the case where we only approximate $\bar{\xb}^{\star}_f$ up to a tolerance level (and not exactly computed), where we can decide on the fly whether such approximation is adequate---see Lemma~\ref{le:choice_of_params} below.

The above observations lead to the following lemma: given a point $\xb$, we bound $\Vert\xb - \xbar^{\star}_t\Vert_{\xbar^{\star}_t}$ by the distance $\Vert\xb - \xb^{\star}_t\Vert_{\xb^{\star}_t}$, using the bound \eqref{eq:Delta_star_est}.  

\begin{lemma}\label{le:Delta_star_est2}
Consider the configuration in Corollary \ref{co:Delta_star_est}, such that $\bar{m}_0 < \tfrac{1}{2}$.
Let $\lambda_t(x) := \Vert x - x_t^{\star}\Vert_{x^{\star}_t}$ and $\bar{\lambda}_t(x) := \Vert x - \bar{x}^{\star}_t\Vert_{\bar{x}^{\star}_t}$, for any $\xb\in\intx{\Xc}$. Then, the following connection between $\lambda_t(\xb)$ and $\bar{\lambda}_t(\xb)$ holds:
\begin{equation}\label{eq:Delta_star_est2}
\bar{\lambda}_t(x) \leq \frac{\lambda_t(x)}{1-\bar{\Delta}_t} + \bar{\Delta}_t 
\leq \frac{(1 - \bar{m}_0)\lambda_t(x)}{1 - 2\bar{m}_0} + \frac{\bar{m}_0}{1-\bar{m}_0}.
\end{equation}
\end{lemma}

\textit{Proof.}
By definition of the local norm $\bar{\lambda}_t(\xb)$, we have
\begin{align*}
\bar{\lambda}_t(\xb) &= \iprods{\nabla^2{f}(\xbar^{\star}_t)(\xb - \xbar^{\star}_t),  \xb - \xbar^{\star}_t}^{1/2} \nonumber\\
&\leq \iprods{\nabla^2{f}(\xbar^{\star}_t)(\xb^{\star}_t - \xbar^{\star}_t),  \xb^{\star}_t - \xbar^{\star}_t}^{1/2} 
+ \iprods{\nabla^2{f}(\xbar^{\star}_t)(\xb - \xb^{\star}_t),  \xb - \xb^{\star}_t}^{1/2} \nonumber\\
&\leq \bar{\Delta}_t + \left(1 - \Vert \xb^{\star}_t - \xbar^{\star}_t \Vert_{\xbar^{\star}_t} \right)^{-1}\iprods{\nabla^2{f}(\xb^{\star}_t)(\xb - \xb^{\star}_t),  \xb - \xb^{\star}_t}^{1/2}\nonumber\\
&= \bar{\Delta}_t + \frac{\lambda_t(\xb)}{1 - \bar{\Delta}_t}.
\end{align*}
Here, in the first inequality, we use the triangle inequality for the weighted norm $\norm{\cdot}_{\nabla^2{f}(\xbar^{\star}_t)}$, while in the second inequality we apply \cite[Theorem 4.1.6]{Nesterov2004}.
The proof is completed when we use \eqref{eq:Delta_star_est} to upper bound the RHS. 
\Eproof

\medskip 
The above lemma indicates that, given a fixed $t > 0$, any approximate solution to \eqref{eq:modified_cvx_prob} (say $\widehat{\xb}_t$)  that is ``good" enough (\emph{i.e.}, the metric $\lambda_t(\widehat{\xb}_t)$ is small), signifies that $\widehat{\xb}_t$ is also ``close'' to the optimal of \eqref{eq:barrier_prob2} (\emph{i.e.}, the metric $\bar{\lambda}_t(\widehat{\xb}_t)$ is bounded by $\lambda_t(\widehat{\xb}_t)$ and, thus, can be controlled).
This fact allows the use of \eqref{eq:modified_cvx_prob}, instead of \eqref{eq:barrier_prob2}, and provides freedom to cleverly select initial parameters $t_0$ and $\eta$ for faster convergence. The next section proposes such an initialization procedure.

\subsection{The choice of initial parameters.}
Here, we describe how we initialize $t_0$ and $\eta$.
Lemma \ref{le:Delta_star_est2} suggests that, for some $\beta \in (0, 1)$, if we can bound $\lambda_t(\cdot) \leq \beta$, then $\bar{\lambda}_t(\cdot)$ is bounded as well as $\bar{\lambda}_t(\cdot) \leq  \frac{(1-\bar{m}_0)\beta}{1 - 2\bar{m}_0} + \frac{\bar{m}_0}{1-\bar{m}_0}$.
This observation leads to the following lemma. 
\begin{lemma}\label{le:bound_of_lbd_0}
Let $\lambda_{t_0}(\xb^0) := \Vert\xb^0 - \xb^{\star}_{t_0}\Vert_{\xb^{\star}_{t_0}}$, where $\xb^{\star}_{t_0}$ is the solution of \eqref{eq:modified_cvx_prob} at $t := t_0$ and $\xb^0$ is an arbitrarily chosen initial point in $\dom{f}$. 
Let $\xi_0\in \partial{g}(\xb^0)$ and, from \eqref{eq:rteta}, $r_{t_0, \eta}(\xb^0) :=  \tfrac{1}{t_0}\left(\cb + \xi_0\right) + \nabla{h}_{\eta}(\xb^0)$. Then, we have
\begin{equation}\label{eq:bound_of_lbd_0}
\lambda_{t_0}(\xb^0) \leq \frac{ 1 - \gamma_{t_0} - \sqrt{ 1 - 6\gamma_{t_0} + \gamma_{t_0}^2} }{2},
\end{equation}
provided that $\gamma_{t_0} :=  \Vert r_{t_0, \eta}(\xb^0) \Vert^{*}_{\xb^0}  \equiv \abs{1-\eta}\Vert \zeta_0\Vert_{\xb^0}^{*} < 3 - 2\sqrt{2}$ for a particular choice of $\eta$.
\end{lemma}

\textit{Proof.}
Since $\xb_{t_0}^{\star}$ is the solution of \eqref{eq:modified_cvx_prob} at $t = t_0$, there exists $\xi_{t_0}^{\star}  \in \partial{g}(\xb^{\star}_{t_0})$ such that: $\nabla h_{\eta}(\xb^{\star}_{t_0}) + \tfrac{1}{t_0} \left(\cb + \xi^{\star}_{t_0}\right) = 0$.
Hence, by the definitions of $r_{t_0, \eta} (\xb^0)$ and $h_\eta(\eta)$, we obtain
\begin{align*}
-\tfrac{1}{t_0} \cdot \xi_{t_0}^* &= \tfrac{1}{t_0} \cdot c + \nabla h_\eta(\xb_{t_0}^*) ~~~~~\Rightarrow~~~~~\tfrac{1}{t_0} \cdot (\xi_0 - \xi_{t_0}^{\star} ) = r_{t_0, \eta}(\xb^0) - \left(\nabla{f}(\xb^0) - \nabla{f}(\xb^{\star}_{t_0})\right).
\end{align*}
By convexity of $g$, we have 
\begin{equation*}
0 \leq \iprods{\xi_0 - \xi_{t_0}^{\star} , \xb^0 - \xb^{\star}_{t_0}} = \iprods{r_{t_0, \eta}(\xb^0) - \left(\nabla{f}(\xb^0) - \nabla{f}(\xb^{\star}_{t_0})\right), \xb^0 - \xb^{\star}_{t_0}}.
\end{equation*} 
This inequality leads to 
$\iprods{ \nabla{f}(\xb^0) - \nabla{f}(\xb^{\star}_{t_0}), \xb^0 - \xb^{\star}_{t_0}} \leq \iprods{r_{t_0, \eta}(\xb^0), \xb^0 - \xb^{\star}_{t_0}}$. 
Using the self-concordance of $f$ in \cite[Theorem 4.1.7]{Nesterov2004} and the Cauchy-Schwarz inequality, 
we can derive
\begin{align*}
\frac{\lambda_{t_0}(\xb^0)^2}{1 + \lambda_{t_0}(\xb^0)} &\leq \iprods{ \nabla{f}(\xb^0) - \nabla{f}(\xb^{\star}_{t_0}), \xb^0 - \xb^{\star}_{t_0}} \leq \iprods{r_{t_0, \eta}(\xb^0), \xb^0 - \xb^{\star}_{t_0}} \leq \Vert r_{t_0, \eta}(\xb^0) \Vert^{*}_{\xb^{\star}_{t_0}}\lambda_{t_0}(\xb^0).
\end{align*} 
Hence, $\tfrac{\lambda_{t_0}(\xb^0)}{1 + \lambda_{t_0}(\xb^0)} \leq \Vert r_{t_0, \eta}(\xb^0) \Vert^{*}_{\xb^{\star}_{t_0}}$. 
Moreover, by \cite[Theorem 4.1.6]{Nesterov2004}, we have $\Vert r_{t_0, \eta}(\xb^0) \Vert^{*}_{\xb^{\star}_{t_0}} \leq \tfrac{\Vert r_{t_0, \eta}(\xb^0) \Vert^{*}_{\xb^0}}{1 - \lambda_{t_0}(\xb^0)}$. 
Combining these two inequalities, we obtain 
\begin{equation*}
\frac{\lambda_{t_0}(\xb^0)(1 - \lambda_{t_0}(\xb^0))}{1 + \lambda_{t_0}(\xb^0)} \leq \Vert r_{t_0, \eta}(\xb^0) \Vert^{*}_{\xb^0} \equiv \gamma_{t_0}.
\end{equation*}
After few elementary calculations, one can see that if $\gamma_{t_0} < 3 - 2\sqrt{2}$, we obtain \eqref{eq:bound_of_lbd_0}, 
which also guarantees its right-hand side of \eqref{eq:bound_of_lbd_0} to be positive.
\Eproof

In plain words, Lemma \ref{le:bound_of_lbd_0} provides a recipe for initial selection of parameters:
Our goal is to choose an initial point $\xb^0$ and the parameters $\eta$ and $t_0$ such that $\lambda_{t_0}(\xb^0) \leq \beta$, for a predefined constant $\beta \in (0, 1)$.
The following lemma $(i)$ provides sufficient conditions---that can be checked \emph{a priori}---on the set of initial points that lead to specific $\eta, t_0$ configurations and make the results in the previous subsection hold; and $(ii)$ suggests that even an approximation of the analytical center $\bar{\xb}^{\star}_f$ is sufficient. 

\begin{lemma}\label{le:choice_of_params}
The initial point $\xb^0\in\dom{f}$ and the parameters $\eta$ and $t_0$ need to satisfy the following relations:
\begin{equation}\label{eq:param_condition}
\eta \in \left[1 - \frac{\beta(1-\beta)}{(1+\beta)\vert\norm{\zeta_0}_{\xb^0}^{*}}, ~1 + \frac{\beta(1-\beta)}{(1+\beta)\vert\norm{\zeta_0}_{\xb^0}^{*}}\right]~~\text{and}~~\eta  < \frac{1-\beta}{(3+\beta)n_{\nu}\Vert \zeta_0\Vert_{\xbar^{\star}_f}^{*}}.
\end{equation}
If we choose $\xb^0\in\dom{f}$ such that
\begin{equation}\label{eq:choice_of_x0}
\Vert\nabla{f}(\xb^0)\Vert_{\xb^0}^{\ast} \leq \kappa < \frac{1}{2}\left(2a_0+1-\sqrt{4a_0^2+1}\right)\leq\frac{1}{2}, ~~\text{with}~~a_0 := \frac{1-\beta}{(3+\beta)n_{\nu}} < \frac{1}{3},
\end{equation}
then one can choose $\eta = 1$ and 
\begin{equation}\label{eq:t0_condition2}
t_0 \geq \frac{\Vert\cb + \xi_0\Vert_{\xb^0}^{\ast}(1-\kappa)(3+\beta)n_{\nu}}{(1-2\kappa)(1-\beta) - \kappa(1-\kappa)(3+\beta)n_{\nu}}.
\end{equation}
In addition, the quantity $\bar{m}_0$ defined in Corollary~\ref{co:Delta_star_est} is bounded by $\bar{m}_0 \leq \hat{m}_0 := \tfrac{(1-\kappa)\left(\kappa + t_0^{-1}\Vert\cb + \xi_0\Vert_{\xb^0}^{\ast}\right)}{(1 - 2\kappa)n_{\nu}}$.
As a special case, if $\xb^0 := \bar{\xb}^{\star}_f$, then $\nabla{f}(\xb^0) = 0$, and we can choose $\eta = 1$ and $t_0 > \tfrac{(1-\beta)\Vert \cb + \xi_0\Vert_{\xb^0}^{*}}{(3+\beta)n_{\nu}}$.
If $\eta$ is chosen from the first interval in  \eqref{eq:param_condition}, then we can take any $t_0 > \tfrac{(1-\beta)\Vert \cb + \xi_0\Vert_{\xb^0}^{*}}{\eta(3+\beta)n_{\nu}}$.
\end{lemma}

\textit{Proof.}
Using \eqref{eq:bound_of_lbd_0}, we observe that in order to satisfy $\lambda_{t_0}(\xb^0) \leq \beta$, it is sufficient to require 
\begin{equation*}
1 - \gamma_{t_0} - \sqrt{1 - 6\gamma_{t_0} + \gamma_{t_0}^2}  \leq 2\beta \quad \Rightarrow \quad \gamma_{t_0} \leq \tfrac{\beta(1-\beta)}{1 + \beta}.
\end{equation*}
Since $\gamma_{t_0} = \vert 1-\eta\vert\norm{\zeta_0}_{\xb^0}^{*}$, the inequality $\gamma_{t_0} \leq \tfrac{\beta(1-\beta)}{1 + \beta}$ further implies 
\begin{equation*}
\vert 1 - \eta\vert \leq  \frac{\beta(1-\beta)}{(1+\beta)\norm{\zeta_0}_{\xb^0}^{*}}.
\end{equation*}
Hence, we obtain the first condition of \eqref{eq:param_condition}.

By our theory and the choice of $m_0$ as in Lemma~\ref{le:Delta_star_est}, it holds $\bar{\lambda}_t(\cdot) \leq  \tfrac{(1-m_0)\beta}{1 - 2m_0} + \tfrac{m_0}{1-m_0} < 1$.
Since $\tfrac{m_0}{1-m_0} \leq \tfrac{m_0}{1 - 2m_0}$, the last condition can be upper bounded as follows:
\begin{equation*}
\frac{(1 - m_0)\beta}{1 - 2m_0} + \frac{m_0}{1 - m_0} \leq \frac{(1 - m_0)\beta + m_0}{1 - 2m_0} < 1.
\end{equation*}
This condition suggests to choose $\eta$ such that $m_0 < \tfrac{1-\beta}{3 + \beta}$.
Let $\xb^0\in\dom{f}$ be an initial point, then, by Corollary \ref{co:Delta_star_est}, we can enforce  $m_0 \leq \bar{m}_0 < \tfrac{1-\beta}{3 + \beta}$. 
This condition leads to
\begin{equation*}
\eta n_{\nu}\Vert \zeta_0\Vert_{\xbar^{\star}_f}^{*}  < \frac{1-\beta}{3+\beta},
\end{equation*}
which implies the second condition of \eqref{eq:param_condition}.

If we take $\eta = 1$, then the second condition of \eqref{eq:param_condition} can be written as $\Vert\zeta_0\Vert_{\bar{\xb}^{\star}_f}^{\ast} \leq \tfrac{1-\beta}{(3+\beta)n_{\nu}}$.
Since $\nabla{f}(\bar{\xb}^{\star}_f) = 0$, similar to the proof of \eqref{eq:lm52_est1}, we can show that $\tfrac{\Vert \xb^0 - \bar{\xb}^{\star}_f\Vert_{\xb^0}}{1 + \Vert \xb^0 - \bar{\xb}^{\star}_f\Vert_{\xb^0}} \leq \Vert \nabla{f}(\xb^0)\Vert_{\xb^0}^{\ast}$. Hence, if $\Vert \nabla{f}(\xb^0)\Vert_{\xb^0}^{\ast} \leq \kappa < \tfrac{1}{2}$, we can show that $\Vert \xb^0 - \bar{\xb}^{\star}_f\Vert_{\xb^0} \leq \frac{\kappa}{1-\kappa} < 1$.
Using this condition, we have 
\begin{equation*}
\Vert \zeta_0\Vert_{\bar{\xb}^{\star}_f}^{\star} \leq \tfrac{\Vert\zeta_0\Vert_{\xb^0}^{\ast}}{1 - \Vert\xb^0-\bar{\xb}^{\star}_f\Vert_{\xb^0}} \leq \tfrac{\Vert\nabla{f}(\xb^0)\Vert_{\xb^0}^{\ast} + t_0^{-1}\Vert\cb + \xi_0\Vert_{\xb^0}^{\ast}}{1 - \Vert\xb^0-\bar{\xb}^{\star}_f\Vert_{\xb^0}} \leq \tfrac{(1-\kappa)\left(\kappa + t_0^{-1}\Vert\cb + \xi_0\Vert_{\xb^0}^{\ast}\right)}{1 - 2\kappa}.
\end{equation*}
The condition $\Vert\zeta_0\Vert_{\bar{\xb}^{\star}_f}^{\ast} \leq \tfrac{1-\beta}{(3+\beta)n_{\nu}}$ is guaranteed if we impose
\begin{equation}\label{eq:lm6_est2}
t_0^{-1}\Vert\cb + \xi_0\Vert_{\xb^0}^{\ast} \leq \frac{(1-2\kappa)(1-\beta)}{(1-\kappa)(3+\beta)n_{\nu}} - \kappa = \frac{(1-2\kappa)(1-\beta) - \kappa(1-\kappa)(3+\beta)n_{\nu}}{(1-\kappa)(3+\beta)n_{\nu}}.
\end{equation}
The right hand side is well-defined if $\tfrac{\kappa(1-\kappa)}{1-2\kappa} < \tfrac{1-\beta}{(3+\beta)n_{\nu}} := a_0$, which is equivalent to $0 < \kappa < \tfrac{1}{2}(2a_0+1 - \sqrt{4a_0^2 + 1}) \leq \tfrac{1}{2}$.
This is exactly \eqref{eq:choice_of_x0}. Next, it is clear that \eqref{eq:lm6_est2} implies \eqref{eq:t0_condition2}.

In particular, if we choose as initial point the \emph{analytical center}, \emph{i.e.}, $\xb^0 = \bar{\xb}^{\star}_f$, then $\nabla{f}(\xb^0) = \nabla{f}(\bar{\xb}^{\star}_f) = 0$. 
By Corollary \ref{co:Delta_star_est}, we can enforce  $m_0 \leq \bar{m}_0 < \tfrac{1-\beta}{3 + \beta}$. This condition leads to
\begin{equation*}
\eta n_{\nu}\Vert \zeta_0\Vert_{\xbar^{\star}_f}^{*} = \frac{\eta n_{\nu}}{t_0}\Vert \cb + \xi_0\Vert_{\xb^0}^{*} < \frac{1-\beta}{3+\beta }.
\end{equation*}
If we take $\eta = 1$, which satisfies the first condition of \eqref{eq:param_condition}, then we can choose $t_0$ such that $t_0 > \tfrac{(1-\beta)\Vert \cb + \xi_0\Vert_{\xb^0}^{*}}{(3+\beta)n_{\nu}}$.
This condition provides a rule to select $t_0$ as indicated the last statement.
\Eproof

\section{The single phase proximal path-following algorithm.}\label{sec:proxPF_alg}
In this section, we present the main ideas of our algorithm. According to the previous section, to solve \eqref{eq:constr_cvx}, one can parameterize the path-following scheme \eqref{eq:barrier_prob2} into \eqref{eq:modified_cvx_prob} and, given proper initialization, solve iteratively \eqref{eq:modified_cvx_prob}---\emph{\emph{i.e.}}, in a path-following fashion, for decreasing values of $t$. 

In the following subsections, we describe schemes to solve \eqref{eq:modified_cvx_prob} up to some accuracy and how the errors, due to approximation, propagate into our theory. 
Based on these ideas, Subsection \ref{subsec:algo} describes the main recursion of our algorithm, along with the update rule for $t$ parameter. Subsection \ref{subsec:stopping} provides a practical stopping criterion procedure, such that an $\varepsilon$-solution is achieved. Subsection \ref{subsec:glue} provides an overview of the algorithm and its theoretical guarantees. 

\subsection{An exact proximal Newton scheme.}\label{subsec:prox_newton}
In our discussions so far, $\xb^{\star}_{t}$ denotes the \emph{exact} solution to \eqref{eq:modified_cvx_prob}, for a given value of paramter $t$. 
Since finding $\xb^{\star}_{t}$ exactly might be infeasible (\emph{e.g.}, there might be no closed-form solution), it is common practice to iteratively solve \eqref{eq:modified_cvx_prob} via first- or second-order Taylor approximations of the smooth part.

In this work, we focus on Newton-type solutions. 
Let $Q(\cdot; y)$ be the second-order Taylor approximation of $h_{\eta}(\cdot)$ around $y$, \emph{i.e.},
\begin{align*}
Q(\xb; y) &:= \iprods{\nabla{h}_{\eta}(y), \xb - y} + \tfrac{1}{2}\iprods{\nabla^2{h}_{\eta}(y)(\xb - y), \xb - y} \\ 
&= \iprods{ \nabla{f}(y) - \eta\zeta_0, \xb - y} + \tfrac{1}{2}\iprods{\nabla^2{f}(y)(\xb - y), \xb - y}.
\end{align*}
Then, $\xb^{\star}_t$ can be obtained by iteratively solving the following subproblem:
\begin{align}\label{eq:cvx_subprob}
\xb_{t}^{+} \longleftarrow  \argmin_{\xb\in\intx{\Xc}}\Big\{ \hat{F}_{t}(\xb; \xb_{t}) := Q(\xb; \xb_t) + \tfrac{1}{t} \cdot G(\xb) \Big\},
\end{align} 
where each iteration is computed with \emph{perfect} accuracy. 
Here, $x_t$ denotes the current estimate and $x_t^{+}$ denotes the next estimate. 
Then, repeating \eqref{eq:cvx_subprob} with \emph{infinite} accuracy, we have:
\begin{align*}
\xb_t^{\infty} \equiv \xb^{\star}_t,
\end{align*}
\emph{i.e.}, we ultimately obtain $\xb^{\star}_t$, for fixed $t > 0$.
We note that, for given point $\xb_t$, we can write the optimality condition of \eqref{eq:cvx_subprob} as follows:
\begin{equation}\label{eq:cvx_subprob_opt}
0 \in \nabla{h_{\eta}}(\xb_t) +  \nabla^2{h_{\eta}}(\xb_t)(\xb^\star_t - \xb_t) + \tfrac{1}{t} \cdot \left(\cb + \partial{g}(\xb^\star_t)\right).
\end{equation}
To solve \eqref{eq:cvx_subprob}, one can use  composite convex quadratic minimization solvers; see, \emph{e.g.}, \cite{Beck2009,Boyd2011,Nesterov2007}. 
The efficiency of such solvers affects the overall analytical complexity and is discussed later in the text (see Subsection~\ref{subsec:overall_comp_complexity}).

\subsection{Inexactness in proximal Newton steps.}\label{subsec:inexact}
In practice, we cannot solve \eqref{eq:cvx_subprob} exactly, but only hope for an approximate solution, up to a given accuracy $\delta > 0$ \cite{kyrillidis2014scalable}. 
The next definition characterizes such inexact solutions.

\begin{definition}\label{de:approx_subsol}
Fix $t > 0$ and let $w$ be an anchor point (in \eqref{eq:cvx_subprob}, $w = \xb_t$). 
Moreover, let $\xb_t^\star$ be the exact solution, obtained by solving \eqref{eq:cvx_subprob} exactly. 
We say that a point $z \in \intx{\Xc}$ is a $\delta$-solution to \eqref{eq:cvx_subprob} if
\begin{equation}\label{eq:approx_subsol}
\Vert z - \xb_t^\star \Vert_{w} \leq \delta,
\end{equation}
for a given tolerance $\delta \geq 0$. We denote this notion by $z :\approx x^{\star}_t$.
\end{definition}

Nevertheless, $\xb^\star_t$ is unknown and, thus, we cannot check the condition \eqref{eq:approx_subsol}.
This condition however holds indirectly, when the following holds \cite{TranDinh2013e}:
\begin{equation}\label{eq:approx_subsol2}
\hat{F}_t(z; w) - \hat{F}_t(\xb^\star_t; w) \leq \tfrac{\delta^2}{2}, 
\end{equation}
where $\hat{F}_t(\cdot; \cdot)$ is defined in \eqref{eq:cvx_subprob}. 
This last condition can be evaluated via several convex optimization algorithms, including first-order methods, \emph{e.g.}, \cite{Beck2009,Nesterov2007}.

We will use these ideas next to define our inexact proximal-Newton path-following scheme.

\subsection{A new, inexact proximal-Newton path-following scheme.}\label{subsec:algo}
Here, we design a new, path-following scheme that operates over the re-parameterized central path in \eqref{eq:modified_cvx_prob}. 
This new algorithm chooses an initial point, as described in Section \ref{sec:central_path}, and selects values for parameter $t$ via a new update rule, that differs from that of \cite{TranDinh2013e}.

At the heart of our approach lies the following recursive scheme: 
\begin{equation}\label{eq:cvx_subprob2}
\left\{\begin{array}{ll}
t_{k+1} &:= t_k + d_k, \vspace{1ex}\\
\xb_{t_{k+1}} &:\approx \textrm{arg}\!\!\!\!\!\displaystyle\min_{\xb\in\intx{\Xc}}\big\{ \hat{F}_{t_{k+1}}(\xb; \xb_{t_k}) := Q(\xb; \xb_{t_k}) + \tfrac{1}{t_{k+1}} \cdot G(\xb) \big\}.
\end{array}\right.
\end{equation}
That is, starting from initial points $t_0$ and $\xb^0 \equiv \xb_{t_0}$, we update the penalty parameter $t_k$ to $t_{k+1}$ via the rule $t_{k+1} := t_k + d_k$, at the $k$-th iteration; see next for details. 
Then, we perform \emph{a single} proximal-Newton iteration, in order to approximate the solution to the minimization problem in \eqref{eq:cvx_subprob2}. 
Observe that, while such a step only approximates the minimizer of \eqref{eq:cvx_subprob2}, by satisfying \eqref{eq:approx_subsol} in our analysis, we can still guarantee convergence close to $\xb^\star$ of \eqref{eq:constr_cvx}. 

\begin{table*}[!t]
\centering
\begin{tabular}{c c c} \toprule
\multicolumn{1}{c}{Notation{\!\!\!}} & \phantom{ab} & \multicolumn{1}{c}{Description} \\
\cmidrule{1-1} \cmidrule{3-3} 
\multicolumn{1}{c}{$\xb^{\star}$} & & \multicolumn{1}{l}{Optimal solution of \eqref{eq:constr_cvx}.} \\ 
\multicolumn{1}{c}{$\xbar_t^{\star}$} & & \multicolumn{1}{l}{Exact solution of \eqref{eq:barrier_prob2} for fixed $t$.} \\ 
\multicolumn{1}{c}{$\xb_t^{\star}$} & & \multicolumn{1}{l}{Exact solution of the parametrized problem \eqref{eq:modified_cvx_prob} for fixed $t$.} \\ 
\multicolumn{1}{c}{$\bar{\xb}_{t_{k+1}}$} & & \multicolumn{1}{l}{Exact solution of \eqref{eq:cvx_subprob2} around $\xb_{t_k}$ for the penalty parameter $t_{k+1}$.{\!\!}} \\ 
\multicolumn{1}{c}{$\xb_{t_{k+1}}$} & & \multicolumn{1}{l}{Inexact solution of \eqref{eq:cvx_subprob2} around $\xb_{t_k}$ for the penalty parameter $t_{k+1}$.{\!\!\!}} \\ 
\multicolumn{1}{c}{$\bar{\lambda}_t(\xb)$} & & \multicolumn{1}{l}{$\bar{\lambda}_t(\xb) := \Vert \xb - \bar{\xb}_t^{\star}\Vert_{\bar{\xb}^{\star}_t}$} \\ 
\multicolumn{1}{c}{$\lambda_t(\xb)$} & & \multicolumn{1}{l}{$\lambda_t(\xb) := \Vert \xb - \xb_t^{\star}\Vert_{\xb^{\star}_t}$.} \\ 
\bottomrule
\end{tabular} 
\vspace{1ex}
\caption{Key notation.}  \label{table:1}
\vskip-0.3cm
\end{table*}

\paragraph{Update rule for penalty parameter.}
We define the local distances
\begin{equation}\label{eq:lambda_Delta}
 \lambda_t(\xb) := \norm{\xb  - \xb^{\star}_t}_{\xb^{\star}_t}, ~~\text{and}~~ \Delta_k := \Vert \xb^{ \star}_{t_{k+1}} - \xb^{\star}_{t_k} \Vert_{\xb^{\star}_{t_{k+1}}},
\end{equation}
for any $t > 0$ and $\xb\in\intx{\Xc}$.
Before we provide a closed-form formula for $d_k$ in the update rule $t_{k+1} := t_k + d_k$ above, we require the following lemma, which reveals the relationship between $\lambda_{t_{k+1}}(\xb_{t_{k+1}})$ and $\lambda_{t_{k+1}}(\xb_{t_k})$, 
as well as the relationship between  $\lambda_{t_k}(\xb_{t_k})$ and $\Delta_k$. 

\begin{lemma}\label{le:key_est_of_prox_nt}
Given $\xb_{t_k}$ and $t_{k+1}$, let $\xb_{t_{k+1}}$ be an approximation of $\xb^{\star}_{t_{k+1}}$ computed by the inexact proximal-Newton scheme \eqref{eq:cvx_subprob2}.
Let $\lambda_{t_{k+1}}(\xb_{t_k}) = \Vert \xb_{t_k} - \xb^{\star}_{t_{k+1}}\Vert_{\xb^{\star}_{t_{k+1}}}$ and $\lambda_{t_{k+1}}(\xb_{t_{k+1}}) = \Vert\xb_{t_{k+1}} - \xb^{\star}_{t_{k+1}}\Vert_{\xb^{\star}_{t_{k+1}}}$, according to \eqref{eq:lambda_Delta}. 
Let $\delta \geq 0$ be the user-defined approximation parameter, according to \eqref{eq:approx_subsol}, where $w \equiv \xb^\star_{t_{k+1}}$ at the $(k+1)$-th iteration.
If $\lambda_{t_{k+1}}(\xb_{t_k})\in [0, 1/9]$, then we have
\begin{equation}\label{eq:PN_key_est1}
\lambda_{t_{k+1}}(\xb_{t_{k+1}}) \leq \tfrac{17}{15} \cdot \delta + 5\lambda_{t_{k+1}}(\xb_{t_k})^2.
\end{equation}
The right-hand side of \eqref{eq:PN_key_est1} is nondecreasing w.r.t. $\lambda_{t_{k+1}}(\xb_{t_k})$ and $\delta \geq 0$.
Furthermore, let $\Delta_k$ be defined by \eqref{eq:lambda_Delta}. Then, we have the following estimate:
\begin{equation}\label{eq:PN_key_est2}
\lambda_{t_{k+1}}(\xb_{t_k}) \leq \frac{\lambda_{t_{k}}(\xb_{t_k})}{1 - \Delta_k} + \Delta_k,
\end{equation}
provided that $\Delta_k < 1$.
\end{lemma}

\textit{Proof.}
It is proved in \cite[Theorem 3.3]{TranDinh2013e} that
\begin{equation}\label{eq:estimate1}
\lambda_{t_{k+1}}(\xb_{t_{k+1}}) \leq \frac{\delta}{1- \lambda_{t_{k+1}}(\xb_{t_k})} + \left(\frac{3 - 2\lambda_{t_{k+1}}(\xb_{t_k})}{1 - 4\lambda_{t_{k+1}}(\xb_{t_k}) + 2\lambda_{t_{k+1}}(\xb_{t_k})^2}\right)\lambda_{t_{k+1}}(\xb_{t_k})^2,
\end{equation}
where $\lambda_{t_{k+1}}(\xb_{t_k}) < 1 - \tfrac{1}{\sqrt{2}}$. 
To obtain this result, we require the properties of the proximity operator ${\rm prox}_g(\cdot)$, according to Lemma \ref{lem:proximal_prop}.
Next, we consider the function $m(q) := \tfrac{3-2q}{1 - 4q + 2q^2}$ for $q \in \big[0, 1 - \tfrac{1}{\sqrt{2}}\big)$, in order to model the second term in the right hand side of \eqref{eq:estimate1}. 
We numerically check that 
if $q \in [0, 1/9]$ then $m(q) \leq 5$.
In this case, we also have 
$\tfrac{1}{1 - q} \leq \tfrac{17}{15}$. 
Using these upper bounds into \eqref{eq:estimate1}, we obtain 
$\lambda_{t_{k+1}}(\xb_{t_{k+1}}) \leq \tfrac{17}{15} \cdot \delta + 5\lambda_{t_{k+1}}(\xb_{t_k})^2$, whenever 
$\lambda_{t_{k+1}}(\xb_{t_k}) \in [0, 1/9]$.

For the \eqref{eq:PN_key_est2}, we have
\begin{align*}
\lambda_{t_{k+1}}(\xb_{t_k}) &= \iprods{\nabla^2{f}(\xb^{\star}_{t_{k+1}})(\xb_{t_k} - \xb^{\star}_{t_{k+1}}),  \xb_{t_k} - \xb^{\star}_{t_{k+1}}}^{1/2} \nonumber\\
&\leq \iprods{\nabla^2{f}(\xb^{\star}_{t_{k+1}})(\xb^{\star}_{t_k} - \xb^{\star}_{t_{k+1}}),  \xb^{\star}_{t_k} - \xb^{\star}_{t_{k+1}}}^{1/2} 
+ \iprods{\nabla^2{f}(\xb^{\star}_{t_{k+1}})(\xb_{t_{k}} - \xb^{\star}_{t_{k}}),  \xb_{t_k} - \xb^{\star}_{t_{k}}}^{1/2} \nonumber\\
&\leq \Delta_k + \left(1 - \Vert \xb^{\star}_{t_k} - \xb^{\star}_{t_{k+1}} \Vert_{\xb^{\star}_{t_{k+1}}} \right)^{-1}\iprods{\nabla^2{f}(\xb^{\star}_{t_{k}})(\xb_{t_{k}} - \xb^{\star}_{t_{k}}),  \xb_{t_k} - \xb^{\star}_{t_{k}}}^{1/2}\nonumber\\
&= \Delta_k + \frac{\lambda_{t_k}(\xb_{t_k})}{1 - \Delta_k}.
\end{align*}
Here, in the first inequality, we use the triangle inequality for the weighted norm $\norm{\cdot}_{\nabla^2{f}(\xb^{\star}_{t_{k+1}})}$, while in the second inequality we apply \cite[Theorem 4.1.6]{Nesterov2004}.
\hfill $\blacksquare$

In words, \eqref{eq:PN_key_est1} states that the quadratic convergence rate of consecutive inexact proximal-Newton steps in \eqref{eq:cvx_subprob2} is preserved per iteration. 
Moreover, \eqref{eq:PN_key_est1} describes how the approximation parameter $\delta$ accumulates as the iteration counter increases (\emph{i.e.}, it introduces an additive error term). 

The next lemma shows how we can bound $\Delta_k$ based on the update rule $t_{k+1} = t_k + d_k$ for $d_k \neq 0$. 
This lemma also provides a rule for $d_k$ selection.

\begin{lemma}\label{le:update_t}
Given constant $c_{\beta} >0$, let $\sigma_{\beta} := \tfrac{c_{\beta}}{(1 + c_{\beta})\sqrt{\nu}}$. 
Then, $\Delta_k$ defined in \eqref{eq:lambda_Delta} satisfies
\begin{equation}\label{eq:bound_dt}
\frac{\Delta_k}{1 + \Delta_k} \leq  \frac{\abs{d_k}}{t_k}\Vert\nabla{f}(\xb_{t_{k+1}}^{\ast})\Vert_{\xb_{t_{k+1}}^{\ast}}^{*} \leq \frac{\abs{d_k}\sqrt{\nu}}{t_k}.
\end{equation} 
Moreover, if we choose $d_k := - \sigma_{\beta} t_k $, the $\Delta_k$ is bounded by $\Delta_k \leq c_{\beta}$.
\end{lemma}

\textit{Proof.}
Since $\xb^{\star}_{t_k}$ and $\xb^{\star}_{t_{k+1}}$ are the exact solutions of \eqref{eq:modified_cvx_prob} at $t = t_k$ and $t_{k+1}$, respectively, we have
\begin{equation*} 
0 \in \nabla{h}_{\eta}(\xb^{\star}_{t_k}) + \tfrac{1}{t_k} \cdot \partial{G}(\xb^{\star}_{t_k})~~\text{and}~~ 0 \in \nabla{h}_{\eta}(\xb^{\star}_{t_{k+1}}) + \tfrac{1}{t_{k+1}} \cdot \partial{G}(\xb^{\star}_{t_{k+1}}).
\end{equation*}
Hence, there exist $\vb_{t_k}^{\star} \in \partial{G}(\xb^{\star}_{t_k})$ and $\vb^{\star}_{t_{k+1}} \in \partial{G}(\xb^{\star}_{t_{k+1}})$ 
such that $\vb_{t_k}^{\star} = - t_k\nabla{h_{\eta}}(\xb^{\star}_{t_k})$ and $\vb^{\star}_{t_{k+1}} = -t_{k+1}\nabla{h}_{\eta}(\xb^{\star}_{t_{k+1}})$.
Using the convexity of $G$, we have 
\begin{equation}
\iprods{t_{k+1}\nabla{h}_{\eta}(\xb^{\star}_{t_{k+1}}) - t_k\nabla{h}_{\eta}(\xb^{\star}_{t_k}), \xb_{t_{k+1}}^{\star} - \xb^{\star}_{t_k} } = - \iprods{\vb_{t_{k+1}}^{\star} - \vb^{\star}_{t_k}, \xb_{t_{k+1}}^{\star} - \xb^{\star}_{t_k}}  \leq 0. \label{eq:lm33_est1}
\end{equation}
By the definition $\nabla{h}_{\eta}$ and the update rule \eqref{eq:cvx_subprob2} of $t_k$, we have
\begin{align}\label{eq:lm33_est2}
t_{k+1}\nabla{h}_{\eta}(\xb^{\star}_{t_{k+1}}) - t_k\nabla{h}_{\eta}(\xb^{\star}_{t_k})  = t_k[\nabla{f}(\xb_{t_{k+1}}^{\star}) - \nabla{f}(\xb_{t_{k}}^{\star})] + d_k\nabla{f}(\xb_{t_{k+1}}^{\star}).
\end{align}
Combining \eqref{eq:lm33_est1} and \eqref{eq:lm33_est2}, then using \cite[Theorem 4.1.7]{Nesterov2004}, the Cauchy-Schwarz inequality, 
and the definition of $\Delta_k$ in \eqref{eq:lambda_Delta} we obtain
\begin{align}\label{lm33_est3}
\frac{t_k\Delta_k^2}{1 + \Delta_k} &\leq t_k\iprods{\nabla{f}(\xb_{t_{k+1}}^{\star}) - \nabla{f}(\xb_{t_{k}}^{\star}), \xb_{t_{k+1}}^{\star} - \xb^{\star}_{t_k}} \nonumber \\ 
&\leq -d_k\iprods{\nabla{f}(\xb_{t_{k+1}}^{\star}), \xb_{t_{k+1}}^{\star} - \xb^{\star}_{t_k}} \nonumber\\
& \leq \abs{d_k}\Vert\nabla{f}(\xb_{t_{k+1}}^{\star})\Vert_{\xb_{t_{k+1}}^{\star}}^{*} \Delta_k,
\end{align}
which implies the first inequality of \eqref{eq:bound_dt}. 
The second inequality of \eqref{eq:bound_dt} follows from the fact that $\Vert\nabla{f}(\xb_{t_{k+1}}^{\star})\Vert_{\xb_{t_{k+1}}^{\star}}^{*} \leq \sqrt{\nu}$ due to \cite[formula 2.4.2]{Nesterov2004}.
The last statement of this lemma is a direct consequence of \eqref{eq:bound_dt}.
\hfill $\blacksquare$

Based on the above, we describe next the main result of this section: 
Assume that the point $\xb_{t_k}$ is in the $\beta$-neighborhood of the inexact proximal-Newton method \eqref{eq:cvx_subprob2}, \emph{i.e.}, $\lambda_{t_k}(\xb_{t_k}) \leq \beta$ for given 
$\beta \in (0, 1/9]$, according to Lemma \ref{le:key_est_of_prox_nt}. 
The following theorem describes a condition on $\Delta_k$ such that $\lambda_{t_{k+1}}(\xb_{t_{k+1}}) \leq \beta$. This in sequence determines the update rule of $t$ values.
The following theorem summarizes this requirement.

\begin{theorem}\label{th:quadratic_convergence_region}
Let $\{\lambda_{t_k}(\xb_{t_k})\}$ be the sequence generated by the inexact proximal-Newton scheme \eqref{eq:cvx_subprob2}. 
For any 
$\beta \in (0, 1/9]$, if we choose $\delta$ and $\Delta_k$ such that
\begin{equation}\label{eq:cond_delta_Delta}
\delta \leq \frac{\beta}{16} ~~~~\text{and}~~~\Delta_k \leq \frac{1}{2}\left(1 + 0.43\sqrt{\beta} - \sqrt{(1-0.43\sqrt{\beta})^2 + 4\beta}\right), 
\end{equation}
then the condition $\lambda_{t_k}(\xb_{t_k}) \leq \beta$ implies $\lambda_{t_{k+1}}(\xb_{t_{k+1}}) \leq \beta$.
Consequently, the penalty parameter $t_k$ is updated by
\begin{equation}\label{eq:update_t}
t_{k+1} := (1-\sigma_{\beta})t_k = \left(1 - \frac{c_{\beta}}{(1+c_{\beta})\sqrt{\nu}}\right)t_k,
\end{equation}
which guarantees that $\Delta_k$ satisfies the condition \eqref{eq:cond_delta_Delta}, where 
\begin{equation*}
c_{\beta} := \frac{1}{2}\left(1 + 0.43\sqrt{\beta} - \sqrt{(1-0.43\sqrt{\beta})^2 + 4\beta}\right) 
\in (0, ~ 1/23].
\end{equation*}
In addition, 
$c^{\max}_{\beta} := \max\set{ c_{\beta} : \beta\in (0, 1/9]} \approx \tfrac{1}{23}$ when $\beta \approx 0.042231$.
\end{theorem}

\textit{Proof.}  
Under the assumption $\lambda_{t_k}(\xb_{t_k}) \leq \beta$ and the first condition \eqref{eq:cond_delta_Delta} with 
$\delta \leq \frac{\beta}{16}$, we can obtain from \eqref{eq:PN_key_est1} and \eqref{eq:PN_key_est2} that
$\lambda_{t_{k+1}}(\xb_{t_{k+1}}) \leq \frac{2\beta}{27} + 5\big(\tfrac{\beta}{1-\Delta_k} + \Delta_k\big)^2$. 
To guarantee $\lambda_{t_{k+1}}(\xb_{t_{k+1}}) \leq\beta$, we have
\begin{align*}
\frac{2\beta}{27} + 5\left(\frac{\beta}{1-\Delta_k} + \Delta_k\right)^2 \leq \beta ~~\Rightarrow~~\frac{\beta}{1-\Delta_k} + \Delta_k \leq 0.43\sqrt{\beta}.
\end{align*}
The last condition implies 
\begin{equation*}
\Delta_k \leq \frac{1}{2}\left(1 + 0.43\sqrt{\beta} - \sqrt{(1-0.43\sqrt{\beta})^2 + 4\beta}\right).
\end{equation*}
This is the second condition of \eqref{eq:cond_delta_Delta}, provided that 
$\beta \in (0, 1/9]$.
The second statement of this theorem follows from \eqref{eq:cond_delta_Delta} and Lemma \ref{le:update_t}, while the last statement is computed numerically.
\Eproof

\subsection{Stopping criterion.}{\label{subsec:stopping}}
We require a stopping criterion that guarantees an $\varepsilon$-solution for \eqref{eq:constr_cvx} according to Definition \ref{de:approx_sol}. 
First, let us define the following quantities
\begin{equation}\label{eq:gamma_01}
\gamma_1 := \frac{(1-\hat{m}_0)\beta}{1-2\hat{m}_0} + \frac{\hat{m}_0}{1-\hat{m}_0}~~~\text{and}~~~\hat{\gamma}_1 := \frac{0.43\sqrt{\beta}(1-\hat{m}_0)}{1-2\hat{m}_0} +  \frac{\hat{m}_0}{1-\hat{m}_0},
\end{equation}
where $\hat{m}_0 := \tfrac{(1-\kappa)\left(\kappa + t_0^{-1}\Vert\cb + \xi_0\Vert_{\xb^0}^{\ast}\right)}{(1 - 2\kappa)n_{\nu}}$, as defined in Lemma~\ref{le:choice_of_params}, and 
$\beta \in (0, 1/9]$.
By the proof of Lemma~\ref{le:choice_of_params}, one can show that $\gamma_1 < 1$ and $\hat{\gamma}_0 < 1$.
Next, we present the following Lemma that provides a stopping criterion for our algorithm below; the proof is provided in the Appendix \ref{apdx:le:opt_cond}.

\begin{lemma}\label{le:opt_cond}
Let $\set{\xb_{t_k}}$ be the sequence generated by \eqref{eq:cvx_subprob2}. 
Then, it holds $\forall k $ that $\set{\xb_{t_k}}\subset \Xc$ and 
\begin{equation}\label{eq:g_dist}
0 \leq G(\xb_{t_{k+1}}) - G^{\star} \leq t_{k+1} \cdot \psi_{\beta}(\nu),
\end{equation}
 where 
\begin{equation}\label{eq:psi_fun}
\psi_{\beta}(\nu) := \nu + \sqrt{\nu}\frac{\gamma_1}{1-\hat{\gamma}_0} + \frac{\hat{\gamma}_0}{(1-\hat{\gamma}_0)^2}\left(\hat{\gamma}_0 + \gamma_1 + \delta\right) + \frac{\delta^2}{2} + \hat{m}_0\gamma_1.
\end{equation}
I.e., for a given tolerance $\varepsilon > 0$, if $t_{k+1} \cdot \psi_{\beta}(\nu) \leq \varepsilon$, then $\xb_{t_{k+1}}$ is an $\varepsilon$-solution of \eqref{eq:constr_cvx}, according to Definition~\ref{de:approx_sol}.
\end{lemma} 

\subsection{Overview of our scheme.}{\label{subsec:glue}}
\begin{algorithm}[!ht]\caption{Single-phase, proximal path-following scheme}\label{alg:pathfollowing2}
\begin{algorithmic}[1]
     \vspace{0.5ex}    
    \STATE{\bfseries Input:} Tolerance $\varepsilon > 0$. 
         \vspace{0.5ex}    
   \STATE {\bfseries Initialization:} 
    \STATE\label{alg:step0}{\hspace{2.5ex}}Choose $\beta \in (0, 1/9]$. Set $a_0 : = \tfrac{(1-\beta)}{(3+\beta)n_{\nu}}$ and $\delta := \frac{\beta}{16}$.
    \vspace{0.5ex}
    \STATE\label{alg:step1}{\hspace{2.5ex}}Find a point $\xb^0\in\dom{f}$ such that $\Vert\nabla{f}(\xb^0)\Vert_{\xb^0}^{\ast} \leq \kappa$ for some $\kappa \in \big(0, \tfrac{1}{2}(2a_0\!+\!1 \!-\! \sqrt{4a_0^2+1})\big)$.
     \vspace{0.5ex}
    \STATE\label{alg:step2}{\hspace{2.5ex}}Compute a subgradient $\xi_0 \in \partial{g}(\xb^0)$ and compute $c_0 := \Vert c + \xi_0\Vert_{\xb^0}^{*}$. 
     \vspace{0.5ex}
    \STATE\label{alg:step3}{\hspace{2.5ex}}Set $\eta := 1$, and choose $t_0 > \tfrac{(1-\beta)c_0}{(3+\beta)n_{\nu}}$.
     \vspace{0.5ex}
    \STATE\label{alg:step4}{\hspace{2.5ex}}Compute  $\hat{m}_0$ from Lemma \ref{le:choice_of_params}, $\hat{\gamma}_0$ and $\gamma_1$ from \eqref{eq:gamma_01}, and $\psi_{\beta}(\nu)$ from \eqref{eq:psi_fun} of Lemma~\ref{le:opt_cond}.
     \vspace{0.5ex}    
    \STATE\label{alg:step6}{\hspace{2.5ex}}Set $c_{\beta} := \tfrac{1}{2}\big[1 + 0.43\sqrt{\beta} - \sqrt{(1-0.43\sqrt{\beta})^2 + 4\beta}\big]$ and $\sigma_{\beta} := \tfrac{c_{\beta}}{(1 + c_{\beta})\sqrt{\nu}}$.
         \vspace{1ex}
    \FOR{$k := 0$ {\bfseries to} $k_{\max}$}
     \vspace{0.5ex}    
    	\STATE\label{alg:step7} If $t_k \psi_{\beta}(\nu) \leq \varepsilon$, then TERMINATE.
     \vspace{0.5ex}    	
	\STATE\label{alg:step8} Update $t_{k+1} := (1 - \sigma_{\beta})t_k$.
     \vspace{0.5ex}    	
	\STATE\label{alg:step9} Perform the inexact full-step proximal-Newton iteration by solving
	\begin{equation*}
	\xb_{t_{k+1}} :\approx \mathrm{arg}{\!}\displaystyle\min_{\xb }\big\{ \hat{F}_{t_{k+1}}(\xb; \xb_{t_k}) := Q(\xb; \xb_{t_k}) + \tfrac{1}{t_{k+1}} \cdot G(\xb) \big\}, \quad \quad \text{up to a given accuracy $\delta$.}
	\end{equation*}
   \ENDFOR
\end{algorithmic}
\end{algorithm}
We summarize the proposed scheme  in Algorithm \ref{alg:pathfollowing2}.  
It is clear that the computational bottleneck lies in Step~\ref{alg:step9}, where we approximately solve a strongly convex quadratic composite subproblem. We comment on this step and its solution in Subsection \ref{subsec:overall_comp_complexity}.
The following theorem summarizes the worst-case iteration-complexity of Algorithm \ref{alg:pathfollowing2}. 

\begin{theorem}\label{th:worst_case_complexity2}
Let $\{(\xb_{t_k}, t_k)\}$ be the sequence generated by Algorithm \ref{alg:pathfollowing2}. 
Then, the total number of iterations required to reach an $\varepsilon$-solution $\xb$ of \eqref{eq:constr_cvx} does not exceed 
\begin{equation}\label{eq:worst_case_bound}
k_{\max} := \left\lfloor \frac{\log\left( \tfrac{\psi_{\beta}(\nu)}{t_0\varepsilon}\right)}{-\log(1 - \sigma_{\beta})}\right\rfloor + 1.
\end{equation}
Thus, the worst-case iteration-complexity of Algorithm \ref{alg:pathfollowing2} is $\mathcal{O}\left(\sqrt{\nu}\log\left(\tfrac{\nu}{t_0\varepsilon}\right)\right)$.
\end{theorem}

\textit{Proof. }
From \eqref{eq:update_t}, we can see that $t_k = (1 - \sigma_{\beta})^kt_0$.  
Hence, to obtain $0 \leq G(\xb_{t_{k+1}}) - G^{\star} \leq \varepsilon$, using \eqref{eq:opt_cond2}, we require 
\begin{equation*}
t_k \geq \frac{\psi_{\beta}(\nu)}{\varepsilon}, ~~\text{or}~~~k \geq \frac{\log\big( \tfrac{\psi_{\beta}(\nu)}{t_0\varepsilon}\big)}{-\log(1 - \sigma_{\beta})}.
\end{equation*}
By rounding up this estimate, we obtain $k_{\max}$ as in \eqref{eq:worst_case_bound}.
We note that $-\log(1 - \sigma_{\beta}) = \mathcal{O}(1/\sqrt{\nu})$.
In addition, by \eqref{eq:psi_fun}, we have $\psi_{\beta}(\nu) =  \mathcal{O}\left(\nu\right)$.
Hence, the worst-case  iteration-complexity of Algorithm \ref{alg:pathfollowing2} is $\mathcal{O}\left(\sqrt{\nu}\log\big(\tfrac{\nu}{t_0\varepsilon}\big)\right)$.
\Eproof

We note that the worst-case iteration-complexity stated in Theorem \ref{th:worst_case_complexity2} is a global worst-case complexity, 
which is different from the one in \cite{TranDinh2013e}. As already mentioned in the Introduction, in the latter case we require
\begin{align*}
\left \lfloor \frac{F_{t_0}(\xb^0) - F_{t_0}(\xb^\star_{t_0})}{\omega \left((1 - \kappa)\beta \right)} \right \rfloor
\end{align*} iterations in \textsc{Phase I}, for arbitrary selected $t_0$ and $\xb^0$, and $\kappa \in (0,1), \beta \in (0, 0.15]$, $\omega(q) = q - \log(1 + q)$, \emph{i.e.}, \textsc{Phase I} has a sublinear convergence rate to the initial point $\xb^0_{t_0}$.

\begin{figure}[!htbp]
\begin{center}
\includegraphics[scale=1,width=0.8\textwidth]{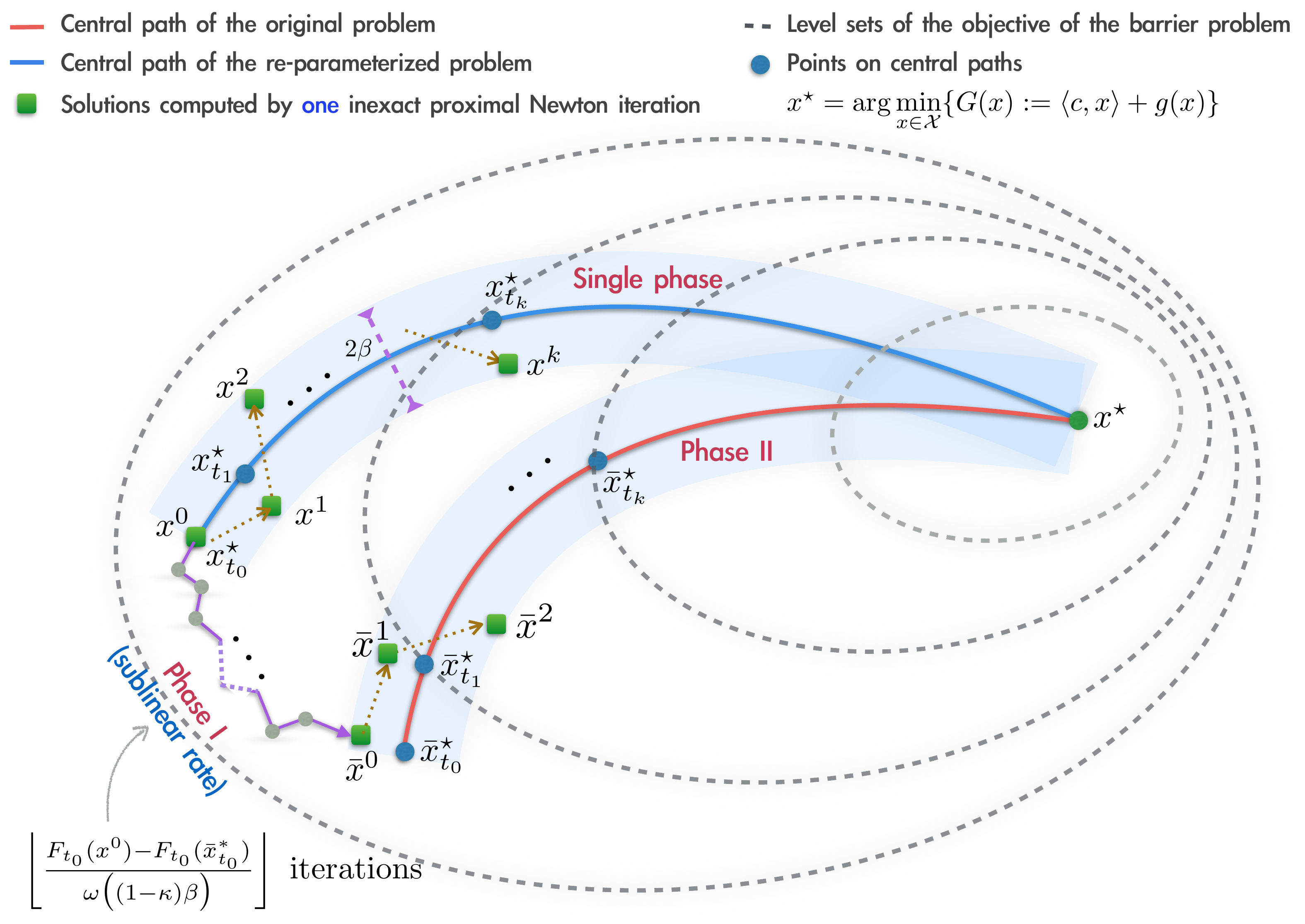}
\end{center}
\caption{Illustration of differences between path-following trajectories followed by single-phase Algorithm \ref{alg:pathfollowing2} and two-phase algorithm in \cite{TranDinh2013e}. In the latter case and given an initial point, say $\xb^ 0 \equiv \xbar_f^\star$, \cite{TranDinh2013e} first performs $\left \lfloor \tfrac{F_{t_0}(\xb^0) - F_{t_0}(\xb^\star_{t_0})}{\omega((1 - \kappa)\beta)} \right \rfloor$ iterations for \textsc{Phase I} to obtain an initial point $\xbar^0$, within the quadratic convergence region of Newton method. Then, the fast convergent \textsc{Phase II} follows the central path (in red color) towards $\xb^\star$. Our algorithm avoids the sublinearly convergent $\textsc{Phase I}$ by properly selecting $t_0, \eta$ and $x^0$, and follows a different central path generated by the solution trajectory of the re-parameterized barrier problem (blue curve).}\label{fig:central_path}
\vspace{-2ex}
\end{figure}

We illustrate the basic idea of our single-phase scheme compared to the two-phase scheme in \cite{TranDinh2013e} in Figure \ref{fig:central_path}. Our method follows different central path generated by the solution trajectory of the re-parameterized barrier problem, where an initial point $x^0$ is immediately available.

\subsection{The exact variant.}
We consider a special case of Algorithm \ref{alg:pathfollowing2}, where the subproblem \eqref{eq:cvx_subprob} at Step 9 can be solved exactly to obtain $\xb_{t_{k+1}}$ such that $\xb_{t_{k+1}} = \xbar^{k+1}_{t_{k+1}}$. 
In this case, we can enlarge the constant $c_{\beta}$ in \eqref{eq:update_t} to obtain a better factor $\sigma_{\beta}$. More precisely, we can use
\begin{equation}\label{eq:exact_c_beta}
\bar{c}_{\beta} := \frac{1}{2}\left[1 + 0.45\sqrt{\beta} - \sqrt{(1-0.45\sqrt{\beta})^2 + 4\beta}\right] > c_{\beta},
\end{equation}
where $\beta \in (0, 0.116764]$.
Hence, we obtain a faster convergence (up to a constant factor) in this case.
For instance, we can numerically check that $\bar{c}_{\beta}^{\max} := \max\set{ \bar{c}_{\beta} : \beta \in (0,  0.116764]} = 0.048186 > c_{\beta}^{\max} = 0.044183$ with respect to $\beta = 0.045864$.

\subsection{Overall computational complexity of Algorithm~\ref{alg:pathfollowing2} in practice.}\label{subsec:overall_comp_complexity}
In this subsection, we present the actual computational complexity of our implementation. 
Key step in the proposed scheme is the subproblem \eqref{eq:cvx_subprob} (or Step~\ref{alg:step9} in Algorithm ~\ref{alg:pathfollowing2}) and how we can efficiently complete it.
Let $\bar{\xb}_{t_{k+1}}$ be the exact solution of this problem, which exists and is unique due to the strong convexity of \eqref{eq:cvx_subprob}. 
Our objective is to solve this problem approximately in order to obtain a $\delta$-approximate solution $\xb_{t_{k+1}}$, according to \eqref{eq:approx_subsol2}.

Observe that the subproblem \eqref{eq:cvx_subprob} contains the second-order Taylor approximation of $h_\eta(\cdot)$; 
in particular, it contains the term $Q(\xb; \xb_{t_k}) := \iprods{\nabla{f}(\xb_{t_k}) - \eta \zeta_0, \xb -\xb_{t_k}} + \tfrac{1}{2}\iprods{\nabla^2{f}(\xb_{t_k})(\xb - \xb_{t_k}), \xb-\xb_{t_k}}$, which is strongly convex. 
The performance of solvers for \eqref{eq:cvx_subprob} depends entirely on the structure of $\nabla^2f(\xb_{t_k})$. 
We consider the following two cases: \vspace{0.2cm}
\begin{itemize}
\item If the solution $\xb^{\star}$ of \eqref{eq:constr_cvx} does not lie on the boundary of $\Xc$ due to the nonlinearity of $g$, then the condition number $\kappa_k :=  \mathrm{cond}\left(\nabla^2f(\xb_{t_k})\right)$ can be bounded by a fixed constant, independent of the dimension $p$ and the solution tolerance $\varepsilon$.
\item If the solution $\xb^{\star}$ of \eqref{eq:constr_cvx} lies on the boundary of $\Xc$, then similar to the proof in \cite{alizadeh1998primal}, the condition number $\kappa_k$ can grow with rate $\mathcal{O}\left(\tfrac{1}{t_k}\right)$.
In this case, the worst-case upper bound on $\kappa_k$ is $\mathcal{O}(\varepsilon^{-1})$ if $\xb_{t_k}$ is an $\varepsilon$-solution since $t_k \leq \varepsilon$. \vspace{0.2cm}
\end{itemize}

We note that, in the second case, the overall computational complexity of Algorithm~\ref{alg:pathfollowing2} could grow, due to the condition number $\kappa_k$. 
However, our algorithm operates directly on the original problem \eqref{eq:constr_cvx} with the dimension of $p$: as we show in the experiments, this is often much better than transforming it into a standard cone problem in a much higher dimensional space.

To approximately solve \eqref{eq:cvx_subprob}, we can apply methods such as accelerated gradient descent, Frank-Wolfe, coordinate descent, and stochastic gradient descent variants.
In the discussion that follows, let us use the well-known fast proximal-gradient method (FISTA) \cite{Beck2009,Nesterov2007}.
Based on the analysis in \cite{Beck2009,Nesterov2007}, in order to guarantee a $\delta$-solution $\xb_{t_{k+1}}$ using FISTA, we require at most
\begin{equation}\label{eq:fista_complexity}
j_k^{\max} := \left\lfloor \sqrt{\kappa_k}\log\left( \frac{\beta(1 + \kappa_k)}{\delta} \right) \right\rfloor + 1 \quad \text{iterations}. 
\end{equation}
The per-iteration complexity is dominated by the following motions: \vspace{0.2cm}
 \begin{itemize}
 \item Lipschitz constant approximation for $Q(\cdot; \cdot)$: this step requires approximating the largest eigenvalue of $\nabla^2f(\xb_{t_k})$, which translates into $\mathcal{O}(p^2)$ computational complexity via a power method.
\item Gradient of $Q(\cdot; \cdot)$: this step requires one matrix-vector multiplication and two vector additions, with computational complexity $\mathcal{O}(p^2)$.
\item Proximal operator $\prox_{t_{k+1}^{-1}g}(\cdot)$: this step requires polynomial time operations as stated in Definition~\ref{de:tractable_proximity}; denote this complexity cost by $\mathcal{T}_{\prox}$. \vspace{0.2cm}
 \end{itemize}
Combining the above, the overall per-iteration complexity for approximately solving \eqref{eq:cvx_subprob} is
\begin{equation}\label{eq:complexity_per_iteration}
\mathcal{O}\left( \sqrt{\kappa_k}\log\left( \frac{\beta(1 + \kappa_k)}{\delta} \right) \max\set{ \mathcal{T}_{\prox}, p^2} \right).
\end{equation}
We note that the linear convergence rate can be also achieved with restarting procedures and without requiring a lower bound of the strong convexity parameter in $Q_k$ \cite{fercoq2016restarting}.

Combining the per-iteration worst-case computational complexity \eqref{eq:complexity_per_iteration} and the worst-case iteration-complexity in Theorem~\ref{th:worst_case_complexity2}, we obtain the overall computational complexity of Algorithm~\ref{alg:pathfollowing2}:
\begin{align*}
\mathcal{O}\left( \sqrt{\max_k \kappa_k }\log\left( \frac{\beta(1 + \max_k \kappa_k)}{\delta} \right) \max\set{ \mathcal{T}_{\prox}, p^2} \cdot \sqrt{\nu}\log\left(\tfrac{\nu}{t_0\varepsilon}\right)\right).
\end{align*}

\section{Numerical experiments.}\label{sec:num_exams}
In this section, we first discuss some implementation aspects of Algorithm \ref{alg:pathfollowing2}: $(i)$ how one can solve efficiently the subproblem in step 9 of Algorithm \ref{alg:pathfollowing2} and, $(ii)$ how we can compute the analytical center $\xbar^\star_f$. In sequence, we illustrate the merits of our approach via three numerical examples, where we compare with state-of-the-art interior-point algorithms.

\paragraph{Inexact proximal-Newton step.} The key step of  Algorithm \ref{alg:pathfollowing2} is  the proximal Newton direction.
This corresponds to solving the following strongly convex quadratic composite problem:
\begin{equation}\label{eq:gen_lasso}
\min_{\db\in\R^p}\set{ q(\db) := \iprods{\hb_k, \db} + \sfrac{1}{2} \cdot \iprods{\Hb_k\db, \db} + g(\xb^k + \db) },
\end{equation}
where $\xb^k, \hb_k\in\R^p$, and $\Hb_k$ is a symmetric positive definite matrix.

There exist efficient first-order and proximal quasi-Newton methods that solve \eqref{eq:gen_lasso}; see \emph{e.g.}, \cite{Beck2009,Becker2011a} for concrete instances of proximal methods, as well as \cite{Tran-Dinh2013a, tran2013proximal} for primal and dual approaches on that matter.
The efficiency of such algorithms strongly depends on the computation of $\prox_g$. 
In addition, since \eqref{eq:gen_lasso} is strongly convex, restart strategies, as in \cite{fercoq2016restarting,Odonoghue2012,Su2014} for first order methods, can lead to faster convergence rates in practice.
When $g$ is absent, \eqref{eq:gen_lasso} reduces to a positive definite linear system $\Hb_k\db = -\hb_k$, which can be efficiently solved by conjugate gradient schemes or Cholesky methods.

\paragraph{The analytical center point.}
To obtain the theoretical complexity bound of Theorem \ref{th:worst_case_complexity2}, we require an approximation of the analytical center $\bar{x}_f^{\star}$ of the barrier function $f$. 
While in most cases approximating $\bar{x}_f^{\star}$ require an iterative procedure, there are several practical cases where $\bar{x}^{\star}_f$ can be computed analytically. 
For example, if $f(\xb) := -\sum_{i=1}^p\log(1-\xb_i^2)$, \emph{i.e.}, $f$ is the barrier of the box set $\Xc := \set{\xb\in\R^p : -1\leq \xb_i \leq 1, ~i=1,\cdots, p}$, then $\bar{\xb}_f^{\star} = \boldsymbol{0} \in \mathbb{R}^p$. 
In the case of $f(X) := -\log\det(X) - \log\det(U - X)$ for the set $\Xc := \set{X\in\Sc_{+}^p : 0 \preceq X \preceq U}$, where $\Sc_{+}^p$ denotes the set of positive semi-definite matrices in $p \times p$ dimensions, we have $\bar{X}_f^{\star} = 0.5U$, where $\bar{X}_f^\star$ denotes the analytical center in a matrix form.
In general cases, $\bar{\xb}_f^{\star}$ can be computed after a few Newton iterations. More details on computation of $\bar{\xb}^{\star}_f$ can be found in \cite{Nesterov2004,Nesterov1994}.

\medskip
Next, we study three numerical examples.
We first compare with the two-phase algorithm in \cite{TranDinh2013e}; 
then, we compare Algorithm \ref{alg:pathfollowing2} with some off-the-shelf interior-point solvers such as SDPT3 \cite{Toh2010}, SeDuMi \cite{Sturm1999} and Mosek \cite{mosek}.

\subsection{The {\rm\textsc{Max-Cut}} problem.}\label{subsec:max_cut}
In this example, we consider the SDP relaxation of the well-known \textsc{Max-Cut} problem. 
In particular, consider the following problem:
\begin{equation}\label{eq:max_cut}
\max_{X}\set{ \tfrac{1}{4}\iprods{L, X} : X \succeq 0, ~\mathrm{diag}(X) = \boldsymbol{e} },
\end{equation}
where $X \in \Sc_{+}^p$ is a positive semi-definite optimization variable, $L$ is the Laplacian matrix of the corresponding underlying graph of the problem, $\mathrm{diag}(X)$ denotes the diagonal of $X$ and $\boldsymbol{e} := (1,1,\cdots, 1)^\top \in \mathbb{R}^p$.
The purpose of this section is to compare Algorithm \ref{alg:pathfollowing2} with the two-phase algorithm in \cite{TranDinh2013e}. We note that in the latter case, the algorithm is also an inexact proximal interior point method, that follows a two-phase procedure.

If we define $c := -\tfrac{1}{4}L$, $g(X) := \delta_{\Omega}(X)$, the indicator of the feasible set $\Omega := \set{ X\in\Sc^p_+ : \mathrm{diag}(X) = \boldsymbol{e}}$, then \eqref{eq:max_cut} can be reformulated into \eqref{eq:constr_cvx}. 
In this case, the proximal operator of $g$ is just the projection onto the affine subspace $\Omega$, which can be computed in a closed form. 
Moreover, \eqref{eq:cvx_subprob} can be solved in a closed from: it requires only one Cholesky decomposition and two matrix-matrix multiplications.

\definecolor{maroon}{cmyk}{1,0.4,0,0}

\begin{table}[!htbp]
\newcommand{\cell}[1]{{\!}#1{\!}}
\newcommand{\cellbf}[1]{{\!}{\color{blue}#1}{\!}}
\rowcolors{2}{white}{maroon!06}
\vspace{-2ex}
\begin{center}
\caption{Summary of results on the small-sized \textsc{Max-Cut} problems. 
Here, $\textrm{Error} := \Vert X^k -X^{\star}_{\textrm{SDPT3}}\Vert_F/\Vert X^{\star}_{\textrm{SDPT3}}\Vert_F$ and $f^{\star}_{\textrm{SDPT3}}$ denotes the objective value obtained by using IPM solver SDPT3 \cite{Toh2010} with high accuracy. For the case of \cite{TranDinh2013e}, the two quantities in Iters column denote the number of iterations required for \textsc{Phase I} and \textsc{Phase II}, respectively.
\texttt{g05\_n.0} is for unweighted graphs with edge probability $0.5$; 
\texttt{pm1s\_100.0} is for a weighted graph with edge weights chosen uniformly from $\{-1, 0, 1\}$ and density $0.1$;
\texttt{wd09\_100.0} is for  a $0.1$ density ten graph with integer edge weights chosen from $[-10,10]$;
\texttt{t2g20\_5555} is for each dimension three two-dimensional toroidal grid graphs with gaussian distributed weights and dimension $20\times 20$;
\texttt{t3g7\_5555} is for each dimension three three-dimensional toroidal grid graphs with gaussian distributed weights and dimension $7\times 7\times 7$.
In these two last problems, the adjacency matrix $\Ab$ is normalized by $\sqrt{\max\abs{A_{ij}}}$.
}\label{tbl:max_cut5}
\begin{small}
\begin{tabular}{ c c c c r r r r c r r r r }
\toprule
\multicolumn{3}{c}{} & & \multicolumn{4}{c}{\cite{TranDinh2013e}} & & \multicolumn{4}{c}{Algorithm \ref{alg:pathfollowing2}} \\ \midrule
Name & $p$ & $f^{\star}_{\textrm{SDPT3}}{~~}$ & & \cell{$f(X)~~~$} &  Error{~~} & Iters & Time[s] & & \cell{$f(X)~~~$} &  Error{~~} & Iters & Time[s]  \\ 
\cmidrule{1-3} \cmidrule{5-8} \cmidrule{10-13}
\cell{g05\_60.0} & \cell{60} & \cell{-59.00} & & \cell{-58.94} &\cell{4.35e-03} & \cell{160/680} & \cell{0.40} & & \cell{-58.94}& \cell{4.35e-03} & \cell{704} & \cellbf{0.32} \\
\cell{g05\_80.0} & \cell{80} & \cell{-80.00} & & \cell{-79.92} &\cell{4.38e-03} & \cell{292/772} & \cell{0.63} & & \cell{-79.92}& \cell{4.39e-03} & \cell{799} & \cellbf{0.48} \\ 
\cell{g05\_100.0} & \cell{100} & \cell{-100.00} & & \cell{-99.90} &\cell{4.41e-03} & \cell{351/877} & \cell{0.94} & & \cell{-99.90}& \cell{4.38e-03} & \cell{910} & \cellbf{0.75} \\ 
\cell{pm1s\_100.0} & \cell{100} & \cell{-52.58} & & \cell{-52.52} &\cell{3.76e-03} & \cell{233/1015} & \cell{1.40} & & \cell{-52.52}& \cell{3.77e-03} & \cell{1042} & \cellbf{0.85} \\  
\cell{w09\_100.0} & \cell{100} & \cell{-80.75} & & \cell{-80.67} &\cell{4.20e-03} & \cell{729/968} & \cell{1.30} & & \cell{-80.67}& \cell{4.21e-03} & \cell{996} & \cellbf{0.87} \\ 
\cell{t3g7\_5555} & \cell{343} & \cell{-20620.30} & & \cell{-20616.76} &\cell{4.45e-03} & \cell{107/32} & \cell{2.23} & & \cell{-20599.78} & \cell{1.99e-03} & \cell{89} & \cellbf{1.30} \\ 
\cell{t2g20\_5555} & \cell{400} & \cell{-31163.19} & & \cell{-31153.93} &\cell{1.33e-02} & \cell{159/99} & \cell{3.41} & & \cell{-31154.04}& \cell{1.24e-02} & \cell{157} & \cellbf{2.21} \\  
\bottomrule
\end{tabular}
\end{small}
\end{center}
\end{table}

We test both algorithms on seven small-sized \textsc{Max-Cut} instances generated by \texttt{Rudy}\footnote{\url{http://biqmac.uni-klu.ac.at/biqmaclib}.}.
We also consider four medium-sized instances from the Gset data set\footnote{\url{http://www.cise.ufl.edu/research/sparse/matrices/Gset/index.html}.}, which were also generated from \texttt{Rudy}.
Both algorithms are tested in a Matlab R2015a environment, running on a MacBook Pro Laptop, 2.6GHz Intel Core i7 with 16GB of memory.  The initial value of $t_0$ is set at $t_0 := 0.025$ for both cases. 
We terminate the execution if $\vert f(X^k) - f^{\star}_{\mathrm{SDPT3}}\vert/\vert f^{\star}_{\mathrm{SDPT3}}\vert \leq 10^{-3}$, where $f(X) := -\mathrm{trace}(LX)$.

\begin{table}[!htbp]
\newcommand{\cell}[1]{{\!}#1{\!}}
\newcommand{\cellbf}[1]{{\!}{\color{blue}#1}{\!}}
\rowcolors{2}{white}{maroon!06}
\begin{center}
\caption{Summary of results on the medium-sized \textsc{Max-Cut} problems. 
Here, $\textrm{Error} := \Vert X^k -X^{\star}_{\textrm{SDPT3}}\Vert_F/\Vert X^{\star}_{\textrm{SDPT3}}\Vert_F$ and $f^{\star}_{\textrm{SDPT3}}$ denotes the objective value obtained by using IPM solver SDPT3 \cite{Toh2010} with high accuracy. For the case of \cite{TranDinh2013e}, the two quantities in Iters column denote the number of iterations required for \textsc{Phase I} and \textsc{Phase II}, respectively.
Each problem Gxx is sparse with \%1 to \%3 upper triangle nonzero, binary entries.
}\label{tbl:max_cut6}
\begin{small}
\begin{tabular}{ c c c c r r r r c r r r r }
\toprule
\multicolumn{3}{c}{} & & \multicolumn{4}{c}{\cite{TranDinh2013e}} & & \multicolumn{4}{c}{Algorithm \ref{alg:pathfollowing2}} \\ \midrule
\cell{Name (Gxx)} & $p$ & $f^{\star}_{\textrm{SDPT3}}{~~}$ & & \cell{$f(X)~~~$} &  Error{~~} & Iters & Time[s] & & \cell{$f(X)~~~$} &  Error{~~} & Iters & Time[s]  \\ 
\cmidrule{1-3} \cmidrule{5-8} \cmidrule{10-13}
\cell{G01} & \cell{800} & \cell{-12080.12} & & \cell{-12080.12} &\cell{1.46e-02} & \cell{149/805} & \cell{104.48} & & \cell{-12080.13}& \cell{1.46e-02} & \cell{569} & \cellbf{62.43} \\ 
\cell{G43} & \cell{1000} & \cell{-7029.29} & & \cell{-7029.30} &\cell{2.03e-02} & \cell{153/1031} & \cell{208.82} & & \cell{-7029.30}& \cell{2.03e-02} & \cell{712} & \cellbf{143.00} \\ 
\cell{G22} & \cell{2000} & \cell{-14116.01} & & \cell{-14116.03} &\cell{3.57e-02} & \cell{215/623} & \cell{805.99} & & \cell{-14116.06}& \cell{3.57e-02} & \cell{561} & \cellbf{741.90} \\ 
\cell{G48} & \cell{3000} & \cell{-5998.57} & & \cell{-5998.57} &\cell{1.38e-02} & \cell{225/2893} & \cell{8487.08} & & \cell{-5998.59}& \cell{1.40e-02} & \cell{1978} & \cellbf{7978.35} \\ \bottomrule
\end{tabular}
\end{small}
\end{center}
\end{table}

The results are provided in Tables \ref{tbl:max_cut5}-\ref{tbl:max_cut6}. Algorithm \ref{alg:pathfollowing2} outperforms \cite{TranDinh2013e} in terms of total computational time, while achieving the same, if not better, solution w.r.t. objective value. 
We observe the following trade-off w.r.t. the algorithm in \cite{TranDinh2013e}: if we increase the initial value of $t_0$ in \cite{TranDinh2013e}, then the number of iterations in Phase I is decreasing, but the number of iterations in Phase II is increasing.
We emphasize that both algorithms use the worst-case update rule without any line-search on the step-size as in off-the-shelf solvers. 

\subsection{The {\rm \textsc{Max}-$k$-\textsc{Cut}} problem.}
Here, we consider the SDP relaxation of the \textsc{Max}-$k$-\textsc{Cut} problem, proposed in \cite[eq. (3)]{frieze1997improved}:
 \begin{equation}\label{eq:max_k_cut}
\max_{X}\set{ \frac{k-1}{2k}\iprods{L, X} : X \succeq 0, ~\mathrm{diag}(X) = \boldsymbol{e}, ~X \geq -\frac{1}{k-1}E_p }.
\end{equation}
Here, $L$ is the Laplacian matrix of the corresponding graph,  $\boldsymbol{e} := (1,1,\cdots, 1)^T$, and $E_p$ is the $p\times p$ all-ones matrix. 
Observe that $X \geq Y$ corresponds to entrywise inequality.
Similarly to \eqref{eq:max_cut}, if we define $c := -\frac{(k-1)}{2k}L$ and $g(X) := \delta_{\Omega}(X)$ with $\Omega := \set{X \in\Sc^p_+ : \mathrm{diag}(X) = \boldsymbol{e}, ~X \geq -\frac{1}{k-1}E_p}$, \eqref{eq:max_k_cut} is a special instance of the class of problems described by \eqref{eq:constr_cvx}. 

We compare Algorithm  \ref{alg:pathfollowing2} with three well-established, off-the-shelf interior-point solvers: SDPT3 \cite{Toh2010}, SeDuMi \cite{Sturm1999}\footnote{Both implementations include Matlab and optimized C-coded parts.}, and Mosek \cite{mosek}\footnote{Available for academic use at \url{https://mosek.com}.}. 
We consider synthetically generated $p$-node graphs, where each edge is generated from a ${\rm Bern}(\sfrac{1}{4}, \sfrac{3}{4})$ probability distribution; we also set $k=4$. 
The parameters of Algorithm  \ref{alg:pathfollowing2} are set as in the previous example, and all algorithms are terminated if  $\vert f(X^k) - f^{\star}\vert/\vert f^{\star}\vert \leq 10^{-5}$, where $f^{\star}$ is the best optimal value produced by three off-the-shelf solvers.
We solve  \eqref{eq:cvx_subprob} with a fast projected gradient method, with adaptive restart and a warm-start strategy \cite{Su2014}: Such configuration requires few iterations to achieve our desired accuracy.

\begin{table}[!htbp]
\begin{center}
\caption{Comparison results on the {\rm \textsc{Max}-$k$-\textsc{Cut}} problem. Here, \texttt{Iters} is the number of iterations; Time[s] is the computational time in second; $f(X) = -\mathrm{trace}(LX)$; \texttt{svars} is the number of slack variables; and \texttt{cnstr} is the number of linear constraints. In addition, we have $p(p+1)/2$ variables in $X$ and one SDP constraint.}\label{tbl:max_k_cut}
\newcommand{\cell}[1]{{\!}#1{\!}}
\newcommand{\cellbf}[1]{\textbf{{\!\!}#1{\!\!}}}
\rowcolors{2}{white}{maroon!06}
\begin{small}
\begin{tabular}{c c c c r r c r r c r r c r r }
\toprule
\multicolumn{1}{c}{\cell{Size}} & \multicolumn{2}{c}{\cell{Lifting}} & & \multicolumn{2}{c}{Algorithm \ref{alg:pathfollowing2}} & & \multicolumn{2}{c}{SeDuMi} & & \multicolumn{2}{c}{SDPT3} & & \multicolumn{2}{c}{Mosek} \\ 
\midrule
$p$ & \cell{\texttt{svars}} & \cell{\texttt{cnstr}} & &  \cell{$f(X)~~~$} &  \cell{Time[s]} & & \cell{$f(X)~~~$} & \cell{Time[s]} & & \cell{$f(X)$~~~} & \cell{Time[s]} & & \cell{$f(X)$~~~} & \cell{Time[s]} \\
\cmidrule{1-3} \cmidrule{5-6} \cmidrule{8-9} \cmidrule{11-12} \cmidrule{14-15}
\cell{50} & \cell{1,225} & \cell{1,275} & & \cell{\mblue{-87.733}} & \cell{7.32} & & \cell{-86.174} & \cell{4.76} & & \cell{-86.160} & \cellbf{2.02} & & \cell{-86.138} & \cell{3.84} \\ 
\cell{75} &  \cell{2,775} & \cell{2,850} & & \cell{\mblue{-166.237}} & \cell{9.80} & & \cell{-166.236} & \cell{55.41} & &  \cell{-166.214} & \cell{10.76} & & \cell{-166.214} & \cellbf{8.91} \\ 
\cell{100} & \cell{4,950} & \cell{5,050}  & & \cell{\mblue{-316.741}} & \cellbf{18.37} & & \cell{-316.746} & \cell{732.16} & & \cell{-316.709} & \cell{48.67} & & \cell{-316.653} & \cell{26.63} \\ 
\cell{150} &   \cell{11,175} & \cell{11,325} & & \cell{\mblue{-654.703}} & \cellbf{73.63} & & \cell{-654.684} & \cell{5,121.34} & & \cell{-654.539} & \cell{484.46} & & \cell{-654.673} & \cell{366.36} \\ 
\cell{200} & \cell{19,900} & \cell{20,100}  & & \cell{\mblue{-1185.784}} & \cellbf{169.08} & & \cell{-1185.783} & \cell{25,521.39} & & \cell{-1185.760} & \cell{2,122.91} & & \cell{-1185.647} & \cell{2,048.95} \\
\bottomrule
\end{tabular}
\end{small}
\end{center}
\vspace{-1ex}
\end{table}

\begin{figure}[!htbp]
\begin{center}
\includegraphics[scale=1,width=1\textwidth]{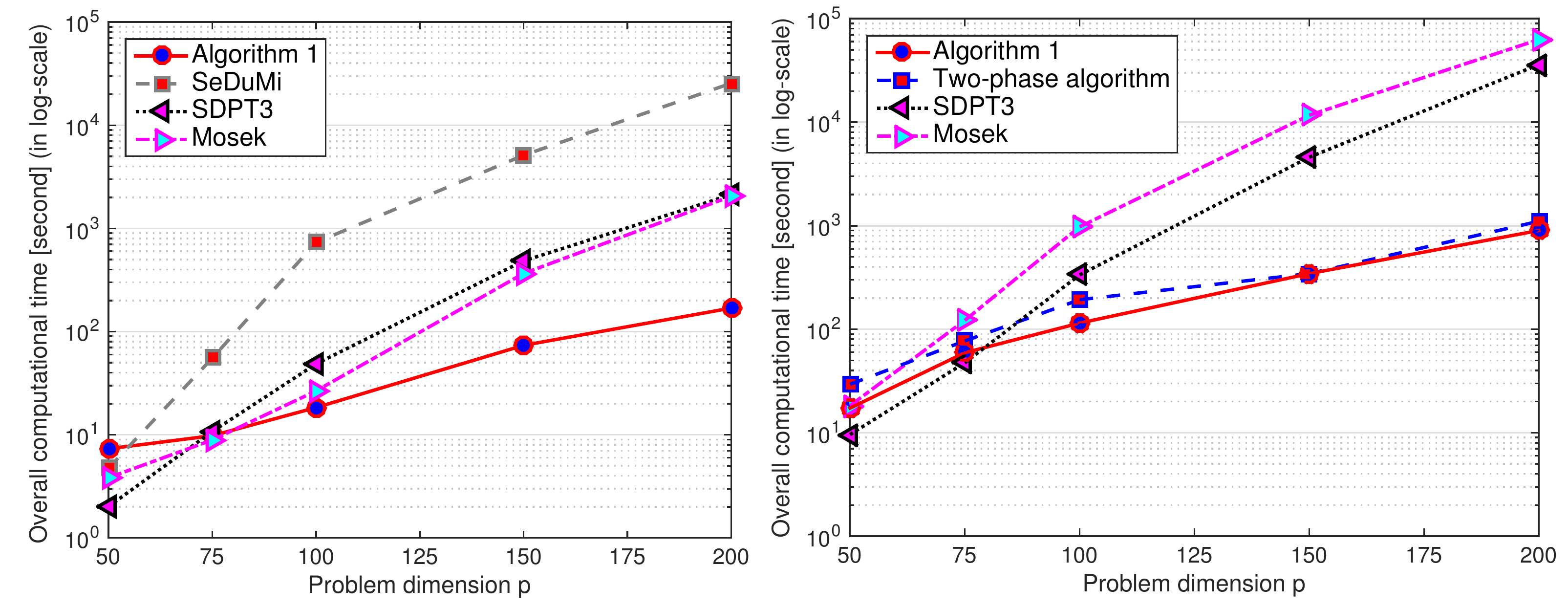}
\end{center}
\caption{Overall execution time, as a function of problem dimension. \textbf{Left panel}: \textsc{Max}-$k$-\textsc{Cut}  problem \eqref{eq:max_k_cut}; \textbf{Right panel}: Clustering problem \eqref{eq:clustering_prob}.}\label{fig:time_scale}
\vspace{-2ex}
\end{figure}

Table \ref{tbl:max_k_cut} contains some experimental results. 
Observe that, if $p$ is small, all algorithms perform well, with the off-the-shelf solvers returning faster a good solution. 
However, when $p$ increases, their computational time significantly increases, as compared to Algorithm \ref{alg:pathfollowing2}. 
One reason that this happens is that standard SDP solvers require $p(p+1)/2$ slack variables and  $p(p+1)/2$ additional linear constraints, in order to process the component-wise inequality constraints. 
Such reformulation of the problem significantly increases variable and constraint size and, hence, lead to slower execution.
In stark contrast, Algorithm \ref{alg:pathfollowing2} handles both linear and inequality constraints by a simple projection, which requires only $p(p+1)/2$ basic operations (flops). 
Figure \ref{fig:time_scale} graphically illustrates the scalability of the four algorithms under comparison, based on the results contained in Table \ref{tbl:max_k_cut}. 

\subsection{Max-norm clustering.}\label{subsec:clustering_prob}
We consider the max-norm clustering task \cite{Jalali2012}, where we seek a clustering matrix $K$ that minimizes the disagreement with a given affinity matrix $\Ab$:
\begin{equation}\label{eq:clustering_prob}
\begin{array}{cl}
\displaystyle\min_{\xb := [L, R, K] \in\mathbb{R}^{p\times 3p}} & \norm{\mathrm{vec}(K - \Ab)}_1\\
\textrm{s.t.} & \mathcal{Q}(\xb) := \begin{bmatrix}L & K \\ K^T & R \end{bmatrix} \succeq 0, ~L_{ii} \leq 1, ~R_{ii}\leq 1, ~i=1,\cdots, p.
\end{array}
\end{equation}
Here, $\mathrm{vec}$ is the vectorization operator of a matrix (\emph{i.e.}, $\mathrm{vec}(\Xb) := (\Xb_1^T, \cdots, \Xb_n^T)^T$, where $\Xb_i$ is the $i$-th column of $\Xb$).
Note that \eqref{eq:clustering_prob} is an SDP convex relaxation to the \emph{correlation clustering} problem; see \cite{Jalali2012} for details.
While \eqref{eq:clustering_prob} comes with rigorous theoretical guarantees and can be formulated as a standard conic program, we need to add $O(p^2)$ slack variables to process the $\ell_1$-norm term and  the linear constraints. 
Moreover, the scaling factors (\emph{e.g.}, the Nesterov-Todd scaling factor regarding the semidefinite cone \cite{Nesterov1997}) can create memory bottlenecks in practice,
by destroying the sparsity of the underlying problem (\emph{e.g.}, by leading to dense KKT matrices in the Newton systems). 

Here, we solve \eqref{eq:clustering_prob} using our path-following scheme. In particular, by defining $\xb := \mathrm{vec}([K, L, R]$, $f(\xb) := -\log\det(\mathcal{Q}(\xb))$ and $g(\xb) := \Vert\mathrm{vec}(K - \Ab)\Vert_1 + \delta_{\mathcal{C}}(\xb)$,  we can transform  \eqref{eq:clustering_prob} into \eqref{eq:constr_cvx}, where $\delta_{\Omega}$ is the indicator function of $\Omega := \set{\xb : L_{ii} \leq 1, ~R_{ii} \leq 1, ~i=1,\cdots, p}$.

We compare the following solvers: Algorithm \ref{alg:pathfollowing2}, the two-phase algorithm in \cite{TranDinh2013e}, SDPT3 and Mosek. 
The initial penalty parameter $t_0$ is set to $t_0 := 0.25$ and the relative tolerance is fixed at $10^{-4}$ for all algorithms.
The data is generated as suggested in \cite{Jalali2012,TranDinh2013e}.
The results of $5$ test problem instances are shown in Table \ref{tbl:clustering_prob} sizes $p$ ranging from $50$ to $200$.\footnote{Since SDPT3 and Mosek cannot run for bigger problems in our personal computer, we restrict to problem sizes up to $p = 200$.}

\begin{table}[!htbp]
\begin{center}
\caption{The performance of Algorithm  \ref{alg:pathfollowing2}, as compared to three methods on the clustering problem \eqref{eq:clustering_prob}. Here, Time[s] is the computational time in second; $f(X) = \Vert\mathrm{vec}(K - A)\Vert_1$, and $s\%$ is the sparsity of $K - A$.}\label{tbl:clustering_prob}
\newcommand{\cell}[1]{{\!\!}#1{\!}}
\rowcolors{2}{white}{maroon!06}
\begin{tabular}{ c c r r r c r r r c r r c r r }
\toprule
\multicolumn{1}{c}{\cell{Size}} & & \multicolumn{3}{c}{Algorithm  \ref{alg:pathfollowing2}} & & \multicolumn{3}{c}{\cite{TranDinh2013e}} & & \multicolumn{2}{c}{SDPT3} & & \multicolumn{2}{c}{Mosek} \\ 
\midrule
$p$ & & \cell{$f(X)~~~$} &  \cell{Time[s]} & s\%{~~}& & \cell{$f(X)~~~$} 	& \cell{Time[s]} & s\%{~~}& & \cell{$f(X)$~~~} & \cell{Time[s]} & & \cell{$f(X)$~~~} & \cell{Time[s]} \\
\cmidrule{1-1} \cmidrule{3-5} \cmidrule{7-9} \cmidrule{11-12} \cmidrule{14-15}
\cell{50}  & & \cell{563.90}  & \cell{17.30} &  \cell{49\%}  &  &  \cell{563.90}  & \cell{29.38} 	& \cell{49.5\%} &  & \cell{563.86} 	& \cell{\mblue{9.60}} &	& \cell{563.92} 		& \cell{18.27} \\ 
\cell{75} & & \cell{1,308.19}    & \cell{59.30}      & \cell{43.8\%}   &  & \cell{1,308.18}    & \cell{77.05} 	& \cell{43.9\%}  & & \cell{1,308.15} 	& \cell{\mblue{47.40}} 	& & \cell{1,308.32} 		& \cell{121.74} \\
\cell{100}   &  & \cell{2,228.62}    & \cell{\mblue{114.59}}  & \cell{35.8\%} & & \cell{2,228.61}    & \cell{192.79} 	&	\cell{35.9\%}        & & \cell{2,228.59} 	& \cell{334.76} 	&	& \cell{2,228.78} 		& \cell{975.10} \\
\cell{150}   &  & \cell{5,328.12}    & \cell{\mblue{344.29}} & \cell{42.4\%} & 	& \cell{5,327.99}    & \cell{344.32}  &	\cell{42.5\%}      & & \cell{5,327.84} 	& \cell{4,584.03} 	&	& \cell{5,328.14} 		& \cell{11,665.52} \\ 
\cell{200}   &  & \cell{9,883.92}    & \cell{\mblue{899.10}}  & \cell{45.8\%} & & \cell{9,883.81}       & \cell{1,102.97} 	& \cell{47.9\%}	& & \cell{9,883.68} 	& \cell{35,974.60} 	&	& \cell{9,884.21} 		& \cell{62,835.42} \\
\bottomrule
\end{tabular}
\end{center}
\vspace{-2ex}
\end{table}

Both SDPT3 and Mosek are approximately $40$ and $60$ times slower than Algorithm~\ref{alg:pathfollowing2} and \cite{TranDinh2013e}, especially when $p > 100$. 
We note that such solvers require $p^2 + 2p$ slack variables, $p^2$ additional second order cone constraints and $2p$ additional linear constraints to reformulate \eqref{eq:clustering_prob} into a standard SDP problem.
Hence, the size of the resulting SDP problem is much larger than of the original one in \eqref{eq:clustering_prob}.
As a concrete example, if $p=200$, then this standard mixed cone problem has a $(200\times 200)$-SDP variable $\Xb$, $40,000$ second order cone variables, $400$ linear variables, $160,400$ linear constraints, and $40,000$ second order cone constraints. 
We also see that Algorithm~\ref{alg:pathfollowing2} is faster than the two-phase algorithm, in terms of total execution time.

We note that SDPT3 gives a slightly better objective value than Algorithm \ref{alg:pathfollowing2}. However, its solution $K$ is fully dense, in contrast to those of Algorithm \ref{alg:pathfollowing2} and \cite{TranDinh2013e}, reducing its interpretation in applications. 
Figure \ref{fig:time_scale} (right) also reveals the scalability of these four algorithms for solving \eqref{eq:clustering_prob}.

\section{Conclusions.}\label{sec:conclusion}
We propose a new path-following framework for a, possibly non-smooth, constrained convex minimization template, which includes linear programming as a special case. 
For our framework, we assume that the constraint set is endowed with a self-concordant barrier and, the non-smooth term has a tractable proximity operator.
Our workhorse is a new re-parameterization of the optimality condition of the convex optimization problem, which allows us to select a different central path towards $\xb^\star$, without relying on the sublinear convergent \textsc{Phase I} of proximal path-following approaches, as in \cite{TranDinh2013e}.

We illustrate that the new scheme retains the same global, worst-case, iteration-complexity with standard approaches \cite{Nesterov2004,Nesterov1994} for smooth convex programming. 
We theoretically show that inexact solutions to sub-problems do not sacrifice the worst-case complexity, when controlled appropriately. 
Finally, we numerically illustrate the effectiveness of our framework on \textsc{Max-Cut} and clustering problems, where the proximal operator play a key role in space efficient optimization.  

\vspace{1ex}
\textbf{Acknowledgments:} 
We would like to thank Yurii Nesterov for useful discussions on the initialization technique used in this work. We also thank Dr. Cong Bang Vu for his overall suggestions on the final version of this manuscript. 
This work was supported in part by NSF grant, no. DMS-16-2044, 
ERC Future Proof, SNF 200021-146750 and SNF CRSII2-147633.

\section{Appendix: The proof of Lemma~\ref{le:opt_cond}.}\label{apdx:le:opt_cond}
For notation, see Table \ref{table:1}.
By \eqref{eq:approx_subsol2} and given $\hat{F}_{t_{k\!+\!1}}(\xb^{k\!+\!1}_{t_{k\!+\!1}}; \xb^k_{t_k}) - \hat{F}_{t_{k\!+\!1}}(\xbar^{k\!+\!1}_{t_{k\!+\!1}}; \xb^k_{t_k}) \leq \tfrac{\delta^2}{2}$, we have
\begin{align}\label{eq:lm35_est1}
G(\xb^{k+1}_{t_{k+1}}) &\leq G(\xbar^{k+1}_{t_{k+1}}) + t_{k+1}Q_k(\xbar^{k+1}_{t_{k+1}};\xb^k_{t_k}) - t_{k+1}Q_{k}(\xb^{k+1}_{t_{k+1}};\xb^k_{t_k}) + \frac{t_{k+1}\delta^2}{2} \nonumber\\
&= G(\xbar^{k+1}_{t_{k+1}}) + t_{k+1}\iprods{\nabla{f}(\xb^k_{t_k}), \xbar^{k+1}_{t_{k+1}} - \xb^{k+1}_{t_{k+1}}} - t_{k+1}\eta\iprods{\zeta_0, \xbar^{k+1}_{t_{k+1}} - \xb^{k+1}_{t_{k+1}}} \nonumber\\
&\quad \quad \quad \quad \quad + \frac{t_{k+1}}{2}\left(\Vert \xbar^{k+1}_{t_{k+1}} - \xb^k_{t_k}\Vert_{\xb^k_{t_k}}^2 - \Vert \xb^{k+1}_{t_{k+1}} - \xb^k_{t_k}\Vert_{\xb^k_{t_k}}^2\right) + \frac{t_{k+1}\delta^2}{2}.
\end{align}
Now, since $\xbar^{k+1}_{t_{k+1}}$ is the exact solution of \eqref{eq:cvx_subprob}, there exists $\bar{\vb}^{k+1} \in \partial{G}(\xbar^{k+1}_{t_{k+1}})$ such that
\begin{equation*}
\bar{\vb}^{k+1}  = -t_{k+1}\big(\nabla{f}(\xb^k_{t_k}) - \eta\zeta_0\big) - t_{k+1}\nabla^2{f}(\xb^k_{t_k})(\xbar^{k+1}_{t_{k+1}} - \xb^k_{t_k}). 
\end{equation*}
Next, using the convexity of $G$, with $\bar{\vb}^{k+1} \in \partial{G}(\xbar^{k+1}_{t_{k+1}})$, we have
\begin{align}\label{eq:lm35_est3}
G(\xbar^{\star}_{t_{k+1}}) - G(\xbar^{k+1}_{t_{k+1}}) &\geq \iprods{\bar{\vb}^{k+1}, \xbar^{\star}_{t_{k+1}} - \xbar^{k+1}_{t_{k+1}}}\nonumber\\
&= -t_{k+1}\iprods{\nabla{f}(\xb^k_{t_k}),  \xbar^{\star}_{t_{k+1}} - \xbar^{k+1}_{t_{k+1}}} + t_{k+1}\eta_0\iprods{\zeta_0,  \xbar^{\star}_{t_{k+1}} - \xbar^{k+1}_{t_{k+1}}} \nonumber\\
&\quad - t_{k+1}\iprods{\nabla^2{f}(\xb^k_{t_k})(\xbar^{k+1}_{t_{k+1}} - \xb^k_{t_k}), \xbar^{\star}_{t_{k+1}} - \xbar^{k+1}_{t_{k+1}}}
\end{align}
Summing up \eqref{eq:lm35_est1} and \eqref{eq:lm35_est3}, and rearranging the terms, we can derive
\begin{align}\label{eq:lm35_est4}
G(\xbar^{\star}_{t_{k+1}}) - G(\xb^{k+1}_{t_{k+1}}) &\geq -t_{k+1}\iprods{\nabla{f}(\xb^k_{t_k}),  \xbar^{\star}_{t_{k+1}} - \xbar^{k+1}_{t_{k+1}}} + t_{k+1}\eta\iprods{\zeta_0,  \xb^{\star}_{t_{k+1}} - \xbar^{k+1}_{t_{k+1}}} \nonumber\\
&\quad - t_{k+1}\iprods{\nabla^2{f}(\xb^k_{t_k})(\xbar^{k+1}_{t_{k+1}} - \xb^k_{t_k}), \xbar^{\star}_{t_{k+1}} - \xbar^{k+1}_{t_{k+1}}}\nonumber\\
&\quad - \frac{t_{k+1}}{2}\left(\Vert \xbar^{k+1}_{t_{k+1}} - \xb^k_{t_k}\Vert_{\xb^k_{t_k}}^2 - \Vert \xb^{k+1}_{t_{k+1}} - \xb^k_{t_k}\Vert_{\xb^k_{t_k}}^2\right) - \frac{t_{k+1}\delta^2}{2}\nonumber\\
&\quad - t_{k+1}\iprods{\nabla{f}(\xb^k_{t_k}), \xbar^{k+1}_{t_{k+1}} - \xb^{k+1}_{t_{k+1}}} + t_{k+1}\eta\iprods{\zeta_0, \xbar^{k+1}_{t_{k+1}} - \xb^{k+1}_{t_{k+1}}} \nonumber\\
&=-t_{k+1}\iprods{\nabla{f}(\xb^k_{t_k}),  \xbar^{\star}_{t_{k+1}} - \xb^{k+1}_{t_{k+1}}} + t_{k+1}\eta\iprods{\zeta_0,  \xbar^{\star}_{t_{k+1}} - \xb^{k+1}_{t_{k+1}}} \nonumber\\
&\quad- t_{k+1}\iprods{\nabla^2{f}(\xb^k_{t_k})(\xbar^{k+1}_{t_{k+1}} - \xb^k_{t_k}), \xbar^{\star}_{t_{k+1}} - \xbar^{k+1}_{t_{k+1}}}\nonumber\\
&\quad - \frac{t_{k+1}}{2}\left(\Vert \xbar^{k+1}_{t_{k+1}} - \xb^k_{t_k}\Vert_{\xb^k_{t_k}}^2 - \Vert \xb^{k+1}_{t_{k+1}} - \xb^k_{t_k}\Vert_{\xb^k_{t_k}}^2\right) - \frac{t_{k+1}\delta^2}{2}.
\end{align}
Now, by the Cauchy-Schwarz inequality, we can further estimate \eqref{eq:lm35_est4} as
\begin{align}\label{eq:lm35_est5}
G(\xbar^{\star}_{t_{k+1}}) - G(\xb^{k+1}_{t_{k+1}}) & \geq -t_{k+1}\Vert\nabla{f}(\xb^k_{t_k})\Vert_{\xbar^{\star}_{t_{k+1}}}^{*} \Vert\xbar^{\star}_{t_{k+1}} - \xb^{k+1}_{t_{k+1}}\Vert_{\xbar^{\star}_{t_{k+1}}} \nonumber\\
&\quad - t_{k+1}\abs{\eta}\Vert \zeta_0\Vert_{\xbar^{\star}_{t_{k+1}}}^{*}\Vert \xbar^{\star}_{t_{k+1}} - \xb^{k+1}_{t_{k+1}}\Vert_{\xbar^{\star}_{t_{k+1}}} \nonumber\\
&\quad - t_{k+1}\iprods{\nabla^2{f}(\xb^k_{t_k})(\xbar^{k+1}_{t_{k+1}} - \xb^k_{t_k}), \xbar^{\star}_{t_{k+1}} - \xbar^{k+1}_{t_{k+1}}}\nonumber\\
&\quad - \frac{t_{k+1}}{2}\left(\Vert \xbar^{k+1}_{t_{k+1}} - \xb^k_{t_k}\Vert_{\xb^k_{t_k}}^2 - \Vert \xb^{k+1}_{t_{k+1}} - \xb^k_{t_k}\Vert_{\xb^k_{t_k}}^2\right) - \frac{t_{k+1}\delta^2}{2}.
\end{align}
We consider the term $$\mathcal{T}_{[1]} := \Vert \xbar^{k+1}_{t_{k+1}} - \xb^k_{t_k}\Vert_{\xb^k_{t_k}}^2 - \Vert \xb^{k+1}_{t_{k+1}} - \xb^k_{t_k}\Vert_{\xb^k_{t_k}}^2 + 2\iprods{\nabla^2{f}(\xb^k_{t_k})(\xbar^{k+1}_{t_{k+1}} - \xb^k_{t_k}), \xbar^{\star}_{t_{k+1}} - \xbar^{k+1}_{t_{k+1}}}.$$
Similarly to the proof of \cite[Lemma 5.1]{TranDinh2013e}, we can show that
\begin{align*} 
\mathcal{T}_{[1]} &\leq \frac{2\bar{\lambda}_{t_{k+1}}(\xb^k_{t_k})}{(1 - \bar{\lambda}_{t_{k+1}}(\xb^k_{t_k}))^2}\left( \bar{\lambda}_{t_{k+1}}(\xb^k_{t_k}) + \bar{\lambda}_{t_{k+1}}(\xb^{k+1}_{t_{k+1}}) + \delta\right).
\end{align*}
Next, by using the self-concordance of $f$ and the definition of $\lambda_{t}(\xb)$, we have
\begin{align}\label{eq:lm35_est7}
\left(1 - \bar{\lambda}_{t_{k+1}}(\xb^k_{t_k})\right)^2\nabla^2{f}(\xbar^{\star}_{t_{k+1}}) \preceq \nabla^2{f}(\xb^k_{t_k}) \preceq \left(1 - \bar{\lambda}_{t_{k+1}}(\xb^k_{t_k})\right)^{-2}\nabla^2{f}(\xbar^{\star}_{t_{k+1}}).
\end{align}
On the one hand, using \eqref{eq:lm35_est7} and $\Vert\nabla{f}(\xb^k_{t_k})\Vert_{\xb^k_{t_k}}^{*} \leq \sqrt{\nu}$, we easily get
\begin{align}\label{eq:lm35_est8}
\Vert\nabla{f}(\xb^k_{t_k})\Vert_{\xbar^{\star}_{t_{k+1}}}^{*} \leq (1 - \bar{\lambda}_{t_{k+1}}(\xb^k_{t_k}))^{-1}\Vert\nabla{f}(\xb^k_{t_k})\Vert_{\xb^k_{t_k}}^{*} \leq (1 - \bar{\lambda}_{t_{k+1}}(\xb^k_{t_k}))^{-1}\sqrt{\nu}.
\end{align}
On the other hand, using $\bar{m}_0$ defined in Corollary~\ref{co:Delta_star_est}, we can show that
\begin{align}\label{eq:lm35_est9}
\Vert \zeta_0\Vert_{\xbar^{\star}_{t_{k+1}}}^{*} \overset{\tiny\text{\cite[Corollary 4.1.7]{Nesterov2004}}}{\leq}  n_{\nu}\Vert\zeta_0\Vert_{\xb_f^{\star}}^{*} \equiv \bar{m}_0 \overset{\tiny\text{Lemma~\ref{le:choice_of_params}}}{\leq} \hat{m}_0.
\end{align}
Substituting \eqref{eq:lm35_est7}, \eqref{eq:lm35_est8} and \eqref{eq:lm35_est9} into \eqref{eq:lm35_est5}, we finally obtain
\begin{align*} 
G(\xbar^{\star}_{t_{k+1}}) - G(\xb^{k+1}_{t_{k+1}}) & \geq -t_{k+1}\Bigg(\frac{\sqrt{\nu}\bar{\lambda}_{t_{k+1}}(\xb^{k+1}_{t_{k+1}})}{(1 - \bar{\lambda}_{t_{k+1}}(\xb^k_{t_k}))^2} + \eta \bar{m}_0\bar{\lambda}_{t_{k+1}}(\xb^{k+1}_{t_{k+1}}) \nonumber\\
&\quad + \frac{\bar{\lambda}_{t_{k+1}}(\xb^k_{t_k})}{(1 - \bar{\lambda}_{t_{k+1}}(\xb^k_{t_k}))^2}\left(\bar{\lambda}_{t_{k+1}}(\xb^k_{t_k}) + \bar{\lambda}_{t_{k+1}}(\xb^{k+1}_{t_{k+1}}) + \delta\right) + \frac{\delta^2}{2}\Bigg).
\end{align*}
Now, let us define the following function
\begin{equation}\label{eq:psi_def}
\psi(\nu, m_0, \hat{\lambda}_0, \lambda_1, \delta) := \nu + \sqrt{\nu}\frac{\lambda_1}{1-\hat{\lambda}_0} + \frac{\hat{\lambda}_0}{(1-\hat{\lambda}_0)^2}\left(\hat{\lambda}_0 + \lambda_1 + \delta\right) + \frac{\delta^2}{2} + m_0\lambda_1.
\end{equation}
Using this definition, and combining $G^{\star} - G(\xbar^{\star}_{t_{k+1}}) \geq -\nu t_{k+1}$ and the last estimate, we obtain
\begin{equation}\label{eq:lm35_est13}
0 \leq G(\xb_{t_{k+1}}) - G^{\star} \leq t_{k+1} \cdot \psi(\nu, \eta\hat{m}_0, \bar{\lambda}_{t_{k+1}}(\xb^{k}_{t_{k}}), \bar{\lambda}_{t_{k+1}}(\xb^{k+1}_{t_{k+1}}), \delta).
\end{equation}
Since $\set{(\xb^k_{t_k}, t_k)}$ is generated by \eqref{eq:cvx_subprob2} for  $\beta \in (0, 1/9]$, we have  $\lambda_{t_{k+1}}(\xb_{t_{k+1}}) \leq \beta$ and $\lambda_{t_{k+1}}(\xb^{k}_{t_{k}}) \leq 0.43\sqrt{\beta}$. 
Using Lemma~\ref{le:Delta_star_est2}, we have
\begin{equation*}
\bar{\lambda}_{t_{k+1}}(\xb_{t_{k+1}}) \leq \frac{(1-m_0)\beta}{1-2m_0} + \frac{m_0}{1-m_0}, ~~\text{and}~~\bar{\lambda}_{t_{k+1}}(\xb^{k}_{t_{k}}) \leq \frac{0.43\sqrt{\beta}(1-m_0)}{1-2m_0} +  \frac{m_0}{1-m_0}.
\end{equation*}
From Lemma~\ref{le:choice_of_params}, we see that $m_0 \leq \hat{m}_0 := \tfrac{(1-\kappa)\left(\kappa + t_0^{-1}\Vert\cb + \xi_0\Vert_{\xb^0}^{\ast}\right)}{(1 - 2\kappa)n_{\nu}}$, we can show that
\begin{equation*}
\bar{\lambda}_{t_{k+1}}(\xb_{t_{k+1}}) \leq \frac{(1-\hat{m}_0)\beta}{1-2\hat{m}_0} + \frac{\hat{m}_0}{1-\hat{m}_0} := \gamma_1~~\text{and}~~ \bar{\lambda}_{t_{k+1}}(\xb^{k}_{t_{k}}) \leq \frac{0.43\sqrt{\beta}(1-\hat{m}_0)}{1-2\hat{m}_0} +  \frac{\hat{m}_0}{1-\hat{m}_0} := \hat{\gamma}_0.
\end{equation*}
Since $\eta := 1$,
using the function $\psi$ defined by \eqref{eq:psi_def}, if we denote $\psi_{\beta}(\nu) := \psi(\nu, \hat{m}_0, \hat{\gamma}_0, \gamma_1, \bar{\delta})$, then
\begin{equation*} 
0 \leq G(\xb_{t_{k+1}}) - G^{\star} \leq t_{k+1} \cdot \psi_{\beta}(\nu),
\end{equation*}
which is exactly \eqref{eq:g_dist}.
Using this estimate, if $t_{k+1} \cdot \psi_{\beta}(\nu) \leq \varepsilon$, then we can say that $\xb_{t_{k+1}}$ is an $\varepsilon$-solution of \eqref{eq:constr_cvx} in the sense of Definition~\ref{de:approx_sol}.
\Eproof

\bibliographystyle{plain}

\end{document}